\documentclass[review,fleqn,3p]{elsarticle} 
\usepackage{amsthm}                     
\usepackage{xcolor,lineno}              
\usepackage{amsmath,amssymb}            
\usepackage{graphicx,threeparttable}    
\usepackage{booktabs,multirow,array}    
\usepackage[colorlinks=true]{hyperref}  
\usepackage[linesnumbered,lined,ruled]{algorithm2e}    

\usepackage{pxfonts}
\usepackage[T1]{fontenc}
\newcommand{\bm}[1]{\boldsymbol{#1}}    

\date{March, 12, 2021}
\begin{document}
\setlength{\baselineskip}{14pt}

\begin{frontmatter}

\title{Sparse dynamical system identification with simultaneous structural parameters and initial condition estimation }

\address[ads1]{College of Economics and Management, Nanjing University of Science and Technology, Nanjing, China}

\author[ads1]{Baolei Wei\corref{corauthor}}
    \cortext[corauthor]{Corresponding author}
    \ead{weiblei@yeah.net}

\begin{abstract}
    Sparse Identification of Nonlinear Dynamics (SINDy) has been shown to successfully recover governing equations from data; however, this approach assumes the initial condition to be exactly known in advance and is sensitive to noise. In this work we propose an integral SINDy (ISINDy) method to simultaneously identify model structure and parameters of nonlinear ordinary differential equations (ODEs) from noisy time-series observations. First, the states are estimated via penalized spline smoothing and then substituted into the integral-form numerical discretization solver, leading to a sparse pseudo linear regression. Then, the sequential threshold least squares is performed to extract the fewest active terms from the overdetermined set of candidate features, thereby estimating structural parameters and initial condition simultaneously and meanwhile, making the identified dynamics parsimonious and interpretable. Simulations detail the method's recovery accuracy and robustness to noise. Examples include a logistic equation, Lokta-Volterra system, and Lorenz system.
\end{abstract}


\end{frontmatter}


\section{Introduction}\label{sec:1}

There is a long and fruit history of modeling time-series data with dynamical systems, resulting in many powerful techniques for system identification \cite{lennart1999system,young2012recursive}. In practice, many branches of science and engineering represent dynamical systems as sets of ODEs, in which each equation describes the change of a single variable as a function of the other variables, and the structures of these ODEs are usually determined according to known physical laws (see e.g. Newton's laws, conservation of mass and energy) that provide first principles from which it is possible to derive equations of motion \cite{ahsan2016differential}. However, it may be difficult or impossible for researchers to derive equations, especially in the emerging fields that are lacking of known physical laws. The problem of discovering governing equation from time-series observations is one of the challenges in modern dynamical system research \cite{brunton2019data}.

There are a host of additional techniques for the discovery of governing equations from time-series data, including the state dependent parameter identification \cite{young2012recursive,young2013hypothetico}, symbolic regression \cite{chen2011timeseries,bongard2007automated,schmidt2009distilling,quade2016prediction}, automated inference \cite{daniels2015automated}, empirical dynamic modelling \cite{ye2015equation}, and deep learning \cite{raissi2018deep,qin2019data}; see \cite{wang2016data} for an incomplete review. These approaches can be viewed as special cases of the fourth paradigm: data-intensive scientific discovery \cite{tolle2011fourth}. Inspired by the key observation that most dynamical systems only consist of a few terms, the recent breakthrough has resulted in a compressive sensing paradigm for the structure identification of underlying nonlinear dynamics \cite{wang2011predicting,schaeffer2013sparse}. By using sparse regression to implement compressive sensing, the SINDy framework was proposed \cite{brunton2016discovering}.
This approach has three principal procedures: (1) estimating the states {together} with derivatives from noisy time-series observations; (2) constructing a library of candidate features to approximate vector fields; (3) conducting sparse regression to generate parsimonious dynamics. At present, there exist a series of outputs about the extensions and applications of SINDy, such as the discovery of partial differential equations \cite{rudy2017data}, the lasso-type sparse optimization \cite{rudy2019data}, and the information criterion-based structure selection \cite{mangan2017model}. However, it should be noticed that SINDy has two limitations: (1) one needs to estimate the accurate derivatives of states from noisy time-series observations, but this may be hard in practice \cite{chou2009recent}, and (2) the initial condition is assumed to be noise-free and exactly known in advance when simulating the fitted trajectory, which can not be satisfied in certain cases.

In a broad sense, the aforementioned SINDy framework can be considered as sparse parameter estimation of ODEs. The subproblem (i.e. estimating parameters of ODEs with known model structures under a measurement error framework) has been a hot topic in the statistical society, especially after Ramsay's work \cite{ramsay2007parameter}; see \cite{ramsay2017dynamic} for a comprehensive review. Integral matching \cite{dattner2015optimal} is an effective approach to this subproblem.
{
    The similar idea to integral matching was first introduced to SINDy by \citet{schaeffer2017sparse} and lead to an integral variant of SINDy. In this variant, the derivative estimation of states is bypassed, one advantage over the standard SINDy. However, both this variant and standard SINDy have the drawbacks that they assume the initial condition to be precisely known in advance rather than estimating it from observations, and they also do not take measurement noise into account, which can lead to misspecified model structures. To this end, the present research improves SINDy framework by integrating integral matching as a parameter estimation approach. The principal contributions are as follows:
    \begin{enumerate}
        \item [(i)]
            The penalised spline smoothing is introduced to estimate dynamical states from noisy observations. This provides a flexible data preprocessing approach to identify the underlying dynamical system from noisy observations.
        \item [(ii)]
            The matrix computation tricks are utilised to separate structural parameters from the initial condition, and a modified sequential threshold least squares is proposed to generate simultaneous estimates of sparse structural parameters and initial condition.
    \end{enumerate}
}

The remainder is structured as follows: section \ref{sec:2} formulates the problem and presents the details of integral SINDy; in section \ref{sec:3} the finite sample performance is demonstrated via simulations and conclusions are presented in section \ref{sec:4}.

\section{Integral SINDy}\label{sec:2}

\subsection{Problem formalization}

Consider the structure identification and parameter estimation problem for homogeneous nonlinear ODE system
\begin{equation}\label{eq:0}
    \begin{cases}
        \frac{d}{dt} {\bm{x}}(t) = \bm{g}\left(\bm{x}(t); \bm{\beta}\right),~ t\geq t_1 \\
        {\bm{x}}(t_1)=  \bm{\eta}
    \end{cases}
\end{equation}
where $\bm{x}(t)\in \mathbb{R}^{d}$ is a $d$-dimensional state vector, $\bm{g}\left( \cdot \right)$ is an unknown vector function, $\bm{\beta}\in \mathbb{R}^{p}$ and $\bm{\eta}\in \mathbb{R}^{d}$ denotes the unknown structural parameters and initial condition, respectively.

In practice, the states are always observed with measurement errors, i.e. the observation equation is
\begin{align}\label{eq:02}
    \bm{y}(t_k)=\bm{x}(t_k)+\bm{e}(k),~ k=1,\cdots,n
\end{align}
where $\bm{e}(k)\sim \mathcal{N}(\bm{0}, \bm{\Sigma})$ is an independent Gaussian noise with mean $\bm{0}$ and finite variance matrix $ \bm{\Sigma}$. Then, for simplicity of presentation, the $n\times d$-dimensional time-series observations are collectively denoted as $\left[\bm{y}_i\right]_{n\times d}$, where each column $\bm{y}_i$, $1\leq i \leq d$, corresponds to the time-series observations $\left\{y_{i}(t_k)\right\}_{k=1}^{n}$.

\begin{figure}[!ht]
    \centering
    \includegraphics[scale=0.85, trim = 0 0 0 0, clip = true]{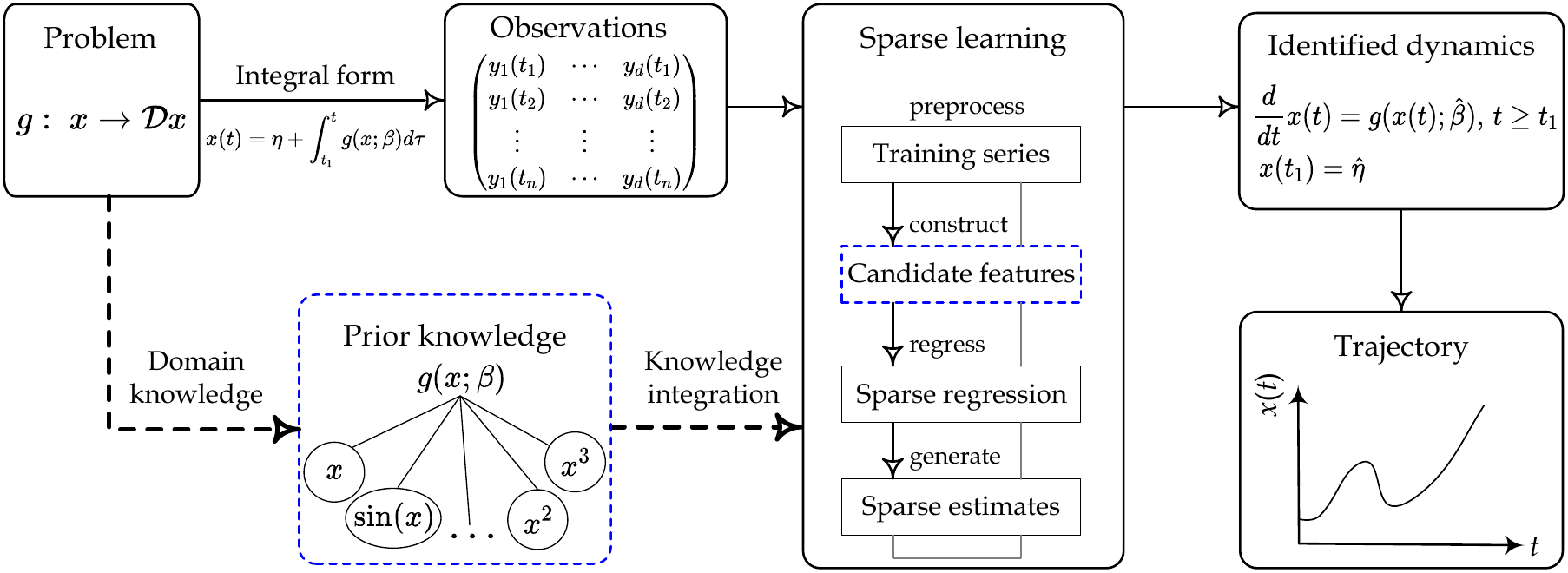} \\
    \caption{Diagrammatic representation of information flow in integral SINDy.}
    \label{fig:1}
\end{figure}

Note that the integral form of state equation \eqref{eq:0} is
\begin{align}\label{eq:01}
    {\bm{x}}(t)=\int_{t_1}^{t} \bm{g}\left(\bm{x}(s); \bm{\beta}\right)ds + \bm{\eta},~ t\geq t_1
\end{align}
which allows for the simultaneous estimation of structural parameter $\bm{\beta}$ and initial condition $\bm{\eta}$. That is, equation \eqref{eq:01} offers an alternative integral SINDy (ISINDy) approach for sparse identification of nonlinear dynamical systems, as shown in Figure \ref{fig:1}.

ISINDy depicts learning from a hybrid information source that consists of data and prior knowledge, also belonging to the emerging field known as informed machine learning \cite{von2021informed,karniadakis2021physics}. The prior knowledge comes from an independent source, such as domain information, empirical experience, and physical laws. The prior knowledge is converted to formal representations --- a library of candidate features, which are used to induce an explicit, likely redundancy expression of vector fields. The following sparse learning consists of three steps: preprocessing the time-series observations to estimate noise-free states, determining the candidate features included in the library, and generating estimates from the sparse regression.

\subsection{ Penalized spline smoothing for states estimation }

Since the observations are contaminated with noise, the penalized spline smoothing \cite{cao2012estimating} is utilized to estimate the states without reference {to} the state equation. For simplicity of presentation, we treat a single-variable case due to that smoothing is applied to each variable in turn. By using a linear combination of basis functions to approximate $x(t)$, it follows that
\begin{equation}\label{eq:05}
    x(t)=\sum_{j=1}^{J} \varphi_j(t) b_j=\bm{\varphi}^\top(t)\bm{b}
\end{equation}
where $\bm{\varphi}(t)=\left[\varphi_1(t),\varphi_2(t),\cdots,\varphi_J(t)\right]^\top$ is the known basis function vector, and $\bm{b}=\left[b_1,b_2,\cdots,b_J\right]^\top$ is the unknown coefficient to be estimated from the observations in equation \eqref{eq:02}. In order to balance the fitness and smoothness of the fitted curve, the objective function contains a penalized sum of squared errors
\begin{equation}\label{eq:06}
    \mathcal{L}(\bm{b})
        =(1-\rho) \sum_{k=1}^{n} \left[y(t_k)-\bm{\varphi}^\top(t_k)\bm{b}\right]^2 +
           \rho \int_{t_1}^{t_n} \left\{ \frac{d^2}{dt^2} \left[ \bm{\varphi}^\top(t)\bm{b} \right] \right\}^2 dt
\end{equation}
where $\rho\in[0, 1)$ is a hyper-parameter; at the right hand side, the first term is the squared error criterion measuring the fitness, and the second term is a smoothing penalty measuring the smoothness. The objective function $\mathcal{L}(\bm{b})$ is a trade-off expression between these two terms: when $\rho=0$, it yields the least-square rough estimates; on the contrary, when $\rho \rightarrow 1$, it yields very smooth estimates and thus nearly linear curve.

Since the basis functions together with their derivatives are known in advance, the penalty term can be reformulated as
\begin{equation}\label{eq:07}
    \int_{t_1}^{t_n} \left\{  \frac{d^2}{dt^2} \left[ \bm{\varphi}^\top(t)\bm{b} \right] \right\}^2 dt
    = \bm{b}^\top \left[\int_{t_1}^{t_n} \frac{d^2}{dt^2}{\bm{\varphi}}(t) \frac{d^2}{dt^2}{\bm{\varphi}}^\top(t) dt\right] \bm{b}
    := \bm{b}^\top \bm{Q} \bm{b}
\end{equation}
where each element of $\bm{Q}$ is calculated as $Q_{i,j}=\int_{t_1}^{t_n} {\frac{d^2}{dt^2}{\varphi}_i(t) \frac{d^2}{dt^2}{\varphi}_j(t)} dt$, and these definite integrals can be approximated using the Simpson's rule.
By substituting equation \eqref{eq:07} into \eqref{eq:06}, partially differentiating the objective function with respect to all parameters, in turn, and setting these derivatives to zero, one obtains the estimates of the coefficients and states:
\begin{equation}\label{eq:08}
    \hat{\bm{b}}=\left[ (1-\rho)\bm{R}^\top\bm{R}+\rho \bm{Q} \right]^{-1}\bm{R}^\top\bm{y}
    \quad \text{and} \quad
    \mathring{\bm{x}}=\bm{R}^\top \hat{\bm{b}}
\end{equation}
where $\bm{R}=[\varphi_j(t_k)]$ is an $n\times J$ matrix, $\bm{y}=\left[ y(t_1), y(t_2), \cdots, y(t_n) \right]^\top$ and $\mathring{\bm{x}}=\left[\mathring{x}(t_1),\mathring{x}(t_2),\cdots,\mathring{x}(t_n) \right]^\top$ are both $n\times 1$ column vectors. To tune the hyper-parameter $\rho$, the cross validation criterion \cite{wahba1990spline} is used, i.e. minimizing the generalized cross validation error
\begin{equation}\label{eq:09}
    \text{GCV}(\rho) = \frac{1}{n}\left\| (\bm{I}-\bm{S}(\rho)) \bm{y} \right\|_2^2
                       / \left[ \frac{1}{n}\mathrm{Trace}(\bm{I}-\bm{S}(\rho)) \right]^2
\end{equation}
where $\bm{S}(\rho)=\bm{R}\left[ (1-\rho)\bm{R}^\top\bm{R}+\rho \bm{Q} \right]^{-1}\bm{R}^\top$ is a $n\times n$ smoothing matrix. It is common to use the grid search strategy to find a reasonable value in the logit scale.

{When implementing the penalised spline smoothing, it is vital to select appropriate basis functions. In practice, the basis functions are often chosen as the cubic B-spline basis functions because of their good properties \cite{cao2012estimating,patrikalakis2002representation}. For example, the property called local compact support, ensures the efficient computation. See Appendix A for details of cubic B-spline basis functions. }

\subsection{Candidate features for vector field representation }

From the prior knowledge that the underling dynamics have only a few terms, the vector field $\bm{g}(\bm{x};\bm{\beta})$ is sparse with respect to the candidate feature space $\bm{\Theta}(\bm{x})$, i.e.
\begin{equation}\label{eq:03}
    \bm{g}\left( \bm{x}, \bm{\beta} \right) = \bm{\Theta}\left( \bm{x} \right) \bm{\Xi}
\end{equation}
where $\bm{\Theta}(\bm{x})=\left[\theta_1(\bm{x})~ \theta_2(\bm{x})~ \cdots~ \theta_m(\bm{x})\right]$ consists of all possible features and $\bm{\Xi}\in\mathbb{R}^{m\times d}$ is a sparse matrix having zeros in each column. Subsequently, without loss of generality, each component of $\bm{g}\left(\bm{x};\bm{\beta}\right)$ can be expressed as a linear combination of the candidate features
\begin{equation}\label{eq:04}
    \left[\bm{g}\left(\bm{x};\bm{\beta}\right)\right]_i
    =\bm{\Theta}(\bm{x})\bm{\xi}_i
    ={\theta}_1(\bm{x})\xi_{1,i}  +  {\theta}_2(\bm{x})\xi_{2,i} + \cdots +  {\theta}_m(\bm{x})\xi_{m,i}
\end{equation}
where $\bm{\xi}_i\in\mathbb{R}^m$, $i=1,\cdots,d$, is the $i$th column of the sparse matrix $\bm{\Xi}$.

By substituting the sparse expression \eqref{eq:03} into state equation \eqref{eq:01}, it follows that
\begin{equation}\label{eq:10}
    \bm{x}(t)=\left(\int_{t_1}^{t} \bm{\Theta}\left( \bm{x}(s) \right)ds\right) \bm{\Xi} + \bm{\eta}, ~
    t\geq t_1
\end{equation}
which, by combining equation \eqref{eq:04}, leads to $d$ separate integral equations
\begin{equation}\label{eq:11}
    x_i(t)=\left( \int_{t_1}^{t} \bm{\Theta}\left( \bm{x}(s) \right)ds\right) \bm{\xi}_i + \eta_i,~
    1\leq i \leq d
\end{equation}
where $\eta_i\in\mathbb{R}$ is the $i$th component of initial condition $\bm{\eta}$. By using the piecewise linear quadrature, the definite integral at each time instant can be approximated as
\begin{equation}\label{eq:12}
    \int_{t_1}^{t_k}\theta_\ell(\bm{x}(t))dt \approx
    \frac{h}{2}\sum_{j=1}^{k-1} \theta_\ell(\bm{x}(t_j)) +
    \frac{h}{2}\sum_{j=2}^{k} \theta_\ell(\bm{x}(t_j)),~
    \ell = 1,2,\cdots,m
\end{equation}
where $h=t_k-t_{k-1}$ is the time interval.

Note that the construction of candidate feature space, $\bm{\Theta}(\bm{x})$, is extremely important to guarantee the accuracy of system identification. Ideally, $\bm{\Theta}(\bm{x})$ includes all terms in the underlying dynamic system; then, $\bm{g}(\bm{x};\bm{\beta})$ can be well represented as a linear combination of the candidate features. But this is not a simple task due to $\bm{\Theta}(\bm{x})$ is not known in advance, even though the prior knowledge can provide some guidance. As a result, it is likely that the number of candidate features, $m$, is large enough to include all possible terms appeared in the vector fields. An alternative strategy is the polynomial and trigonometric functions which has been commonly used in various literature \cite{brunton2019data}.

To illustrate the construction process of $\bm{\Theta}(\bm{x})$, we give an example that has three state variables. Denoting the state variable as $\bm{x}=[x_1~ x_2~ x_3]^\top$, the 3 linear features are $\theta_1(\bm{x})=x_1,~ \theta_2(\bm{x})=x_2,~ \theta_3(\bm{x})=x_3$;
the 6 quadratic features are
$
\theta_4(\bm{x})=x_1x_1,~  \theta_5(\bm{x})=x_1x_2,~  \theta_6(\bm{x})=x_1x_3,~ \theta_7(\bm{x})=x_2x_2,~  \theta_8(\bm{x})=x_2x_3,~  \theta_9(\bm{x})=x_3x_3;
$
and, the 10 cubic features are
$
\theta_{10}(\bm{x})=x_1x_1x_1,~ \theta_{11}(\bm{x})=x_1x_1x_2,~ \theta_{12}(\bm{x})=x_1x_1x_3,~  \theta_{13}(\bm{x})=x_1x_2x_2,~ \theta_{14}(\bm{x})=x_1x_2x_3,~ \theta_{15}(\bm{x})=x_1x_3x_3,~ \theta_{16}(\bm{x})=x_2x_2x_2,~ \theta_{17}(\bm{x})=x_2x_2x_3,~ \theta_{18}(\bm{x})=x_2x_3x_3,~ \theta_{19}(\bm{x})=x_3x_3x_3
$.
Likewise, it is easy to show the trigonometric features; the cosine-based ones are
$
\theta_{20}(\bm{x})=\cos(\omega x_1),~ \theta_{21}(\bm{x})=\cos(\omega x_2),~ \theta_{22}(\bm{x})=\cos(\omega x_3),
$
where $\omega$ is a hyper-parameter. Finally, we obtain a library of candidate features $\bm{\Theta}(\bm{x})=\left[\theta_\ell(\bm{x})\right]$, $\ell=1,\cdots,22$.

Here, if there exists no specific priori knowledge of candidate features, the polynomial functions are suggested because of the following reasons: (1) it is hard to obtain the priori coefficients in the trigonometric functions (see e.g. the hyper-parameter, $\omega$, in $\cos(\omega x_i)$) and (2) in a finite time interval, there always exists a truncated Taylor series to approximate the vector filed globally.

\subsection{Sparse identification by sequential threshold least squares}

In the following context, the candidate features are set to polynomial functions unless otherwise stated, that is, $\theta_\ell(\bm{x})=\prod_{i=1}^{d} x_i^{r_i}$, $1\leq \sum_{i=1}^{d} r_i \leq p$, $r_i\in\mathbb{N}$, and $p$ is the largest degree of polynomials and traditionally set to be less than 3 for the sake of numerical stability.

Denoting the smoothed time series as a matrix form
\[
    \mathring{\bm{X}}=
    \begin{pmatrix}
        | & | &  &  | \\
        \mathring{\bm{x}}_{1} & \mathring{\bm{x}}_{2} & \cdots & \mathring{\bm{x}}_{d} \\
        | & | &  &  |
    \end{pmatrix}
    =\begin{pmatrix}
        \mathring{x}_{1}(t_1) & \mathring{x}_{2}(t_1) & \cdots & \mathring{x}_{d}(t_1) \\
        \mathring{x}_{1}(t_2) & \mathring{x}_{2}(t_2) & \cdots & \mathring{x}_{d}(t_2) \\
        \vdots & \vdots & \cdots & \vdots \\
        \mathring{x}_{1}(t_n) & \mathring{x}_{2}(t_n) & \cdots & \mathring{x}_{d}(t_n)
    \end{pmatrix}
\]
and, by using the element-wise product (denoted by $\odot$), we can obtain the data vectors corresponding to the candidate features. For example, if $\theta_\ell(\bm{x})=x_1x_2x_3$, then   $\bm{\theta}_\ell(\mathring{\bm{X}})=\mathring{\bm{x}}_{1}\odot\mathring{\bm{x}}_{2}\odot\mathring{\bm{x}}_{3}$ and subsequently, the data matrix is $\bm{\Theta}(\mathring{\bm{X}})=\left[\bm{\theta}_1(\mathring{\bm{X}}),~ \cdots,~ \bm{\theta}_m(\mathring{\bm{X}})\right]$ and for ease of notation, denoted by $\mathring{\bm{\Theta}}$.

Substituting data matrices into equation \eqref{eq:11} gives a pseudo linear regression whose estimates can be obtained by minimizing the least-square objective function
\begin{align}\label{eq:13}
    \mathop{\arg\min}_{\bm{\xi}_i,{\eta}_i}
    \left\|
    \bm{S}\mathring{\bm{x}}_i - \bm{C}\mathring{\bm{\Theta}} \bm{\xi}_i - \bm{1} \eta_i
    \right\|_2^2
\end{align}
where $\bm{S}=\left[\bm{0},~ \bm{I}_{n-1}\right]$ selects the last $n-1$ rows, $\bm{C}=\frac{h}{2}\left[ \bm{L}_{n-1},~\bm{0} \right] + \frac{h}{2}\left[\bm{0},~ \bm{L}_{n-1} \right]$ denotes the cumulative operation, $\bm{I}_{n-1}$ is an identity matrix and $\bm{L}_{n-1}$ is a special left triangular with all the entries above the main diagonal are zero and ones elsewhere. It is easy to show the simultaneous least-square estimates of structural parameters and initial condition are
\begin{align}\label{eq:14}
    \begin{pmatrix}
        \hat{\bm{\xi}}_i \\
        \hat{\eta}_i
    \end{pmatrix} =
    \begin{pmatrix}
        \mathring{\bm{\Theta}}^\top \bm{C}^\top\bm{C}\mathring{\bm{\Theta}} &
        \mathring{\bm{\Theta}}^\top \bm{C}^\top\bm{1} \\
        \bm{1}^\top\bm{C}\mathring{\bm{\Theta}} & \bm{1}^\top\bm{1}
    \end{pmatrix}^{-1}
    \begin{pmatrix}
        \mathring{\bm{\Theta}}^\top \bm{C}^\top \bm{S}\mathring{\bm{x}}_i \\
        \bm{1}^\top \bm{S}\mathring{\bm{x}}_i
    \end{pmatrix}
\end{align}
then, according to the inverse formula of block matrix \cite{lu2002inverses}, the inverse matrix can be rewritten as
\begin{align}\label{eq:15}
    \begin{pmatrix}
        \left( \mathring{\bm{\Theta}}^\top \bm{C}^\top\bm{C} \mathring{\bm{\Theta}} \right)^{-1} +
        \frac{1}{\varrho} \left( \mathring{\bm{\Theta}}^\top \bm{C}^\top\bm{C} \mathring{\bm{\Theta}} \right)^{-1} \mathring{\bm{\Theta}}^\top \bm{C}^\top \bm{1} \bm{1}^\top \bm{C} \mathring{\bm{\Theta}} \left( \mathring{\bm{\Theta}}^\top \bm{C}^\top \bm{C}\mathring{\bm{\Theta}} \right)^{-1} &
        -\frac{1}{\varrho} \left( \mathring{\bm{\Theta}}^\top \bm{C}^\top \bm{C}\mathring{\bm{\Theta}} \right)^{-1} \mathring{\bm{\Theta}}^\top \bm{C}^\top \bm{1}  \\
        -\frac{1}{\varrho} \bm{1}^\top \bm{C} \mathring{\bm{\Theta}} \left( \mathring{\bm{\Theta}}^\top \bm{C}^\top \bm{C}\mathring{\bm{\Theta}} \right)^{-1} & \frac{1}{\varrho}
    \end{pmatrix}
\end{align}
where
\[
    \varrho
    = \bm{1}^\top
    \left[
        \bm{I}_{n-1}-
        \bm{C} \mathring{\bm{\Theta}}
        \left( \mathring{\bm{\Theta}}^\top \bm{C}^\top \bm{C}\mathring{\bm{\Theta}} \right)^{-1}
        \mathring{\bm{\Theta}}^\top \bm{C}^\top
    \right] \bm{1}.
\]

By plugging equation \eqref{eq:15} back into \eqref{eq:14}, we can separate the estimates of structural parameters and initial condition, respectively expressed as
\begin{align}\label{eq:16}
    \hat{\bm{\xi}}_i
    =
    \left( \mathring{\bm{\Theta}}^\top \bm{C}^\top\bm{C} \mathring{\bm{\Theta}} \right)^{-1}
    \mathring{\bm{\Theta}}^\top \bm{C}^\top
    \left\{
    \bm{I}_{n-1} -
    \frac{1}{\varrho} \bm{1} \bm{1}^\top
    \left[
        \bm{I}_{n-1} -
        \bm{C} \mathring{\bm{\Theta}}
        \left( \mathring{\bm{\Theta}}^\top \bm{C}^\top\bm{C} \mathring{\bm{\Theta}} \right)^{-1}
        \mathring{\bm{\Theta}}^\top \bm{C}^\top
    \right]
    \right\} \bm{S}\mathring{\bm{x}}_i
\end{align}
and
\begin{align}\label{eq:17}
    \hat{\eta}_i =
    \frac{1}{\varrho} \bm{1}^\top
    \left[
        \bm{I}_{n-1} -
        \bm{C} \mathring{\bm{\Theta}}
        \left( \mathring{\bm{\Theta}}^\top \bm{C}^\top\bm{C} \mathring{\bm{\Theta}} \right)^{-1}
        \mathring{\bm{\Theta}}^\top \bm{C}^\top
    \right]
    \bm{S} \mathring{\bm{x}}_i.
\end{align}

Note that all of the the resultant least-square estimates of structural parameters are not equal to zero and thus cannot select the correct model structure, that is, equation \eqref{eq:16} always results in a full model including all candidate features. In order to perform structure selection, the objective function \eqref{eq:13} is constrained to generate sparse estimates, that is
\begin{align}\label{eq:18}
    \mathop{\arg\min}_{\bm{\xi}_i,{\eta}_i}
        \left\| \bm{S}\mathring{\bm{x}}_i - \bm{C}\mathring{\bm{\Theta}} \bm{\xi}_i - \bm{1}\eta_i \right\|_2^2
    \quad \text{subject to} \quad
        \left\|\bm{\xi}_i\right\|_0 \leq \lambda_i
\end{align}
where $\lambda_i\geq0$ is a hyper-parameter that measures the sparsity of the parameter vector $\bm{\xi}_i$.

By selecting an appropriate hyper-parameter $\lambda_i$, equation \eqref{eq:18} leads to a parsimonious model structure that balances model complexity with descriptive capability. Equation \eqref{eq:18} can be solved by many algorithms, such as the least absolute shrinkage and selection operator (LASSO) \cite{tibshirani2011regression} and least square approximation (LSA) \cite{wang2007unified}. Here, we select the sequential threshold least squares (STLS) \cite{zhang2019convergence} to generate sparse structural estimates, as shown in Algorithm \ref{alg:1}.

\begin{algorithm}[!ht]
    \SetKwInput{Input}{Input}
    \SetKwInOut{Output}{Output}
    \Input{smoothed state matrix $\mathring{\bm{X}}$; candidate data matrix $\mathring{\bm{\Theta}}$; threshold vector $\bm{\lambda}$}

    \For{$i=1,2,\cdots,d$}{
        $\bm{P}^{-1}=\left( \mathring{\bm{\Theta}}^\top \bm{C}^\top\bm{C} \mathring{\bm{\Theta}} \right)^{-1}$;
            \hfill{$\rhd$ initialization based on full model} \\
        ${\varrho} = {n-1} -
                    \bm{1}^\top\bm{C}
                    \mathring{\bm{\Theta}} \bm{P}^{-1} \mathring{\bm{\Theta}}^\top
                    \bm{C}^\top\bm{1} $; \\
        $\hat{\bm{\xi}}_i = \bm{P}^{-1} \mathring{\bm{\Theta}}^\top \bm{C}^\top
        \left\{ \bm{I}_{n-1} - \frac{1}{{\varrho}} \bm{1} \bm{1}^\top
                \left[
                    \bm{I}_{n-1} -
                    \bm{C} \mathring{\bm{\Theta}} \bm{P}^{-1} \mathring{\bm{\Theta}}^\top \bm{C}^\top
                \right]
        \right\} \bm{S}\mathring{\bm{X}}(:,i)$; \\
        \Repeat{convergence}{
        $\text{binds} = (\text{abs}(\hat{\bm{\xi}}_i)>\lambda_j)$;
            \hfill{$\rhd$ threshold of small coefficients} \\
        $\hat{\bm{\xi}}_i(^\sim\text{binds}) = 0$; \\
        $\mathring{\bm{\Omega}} = \mathring{\bm{\Theta}}(:,\text{binds})$;
            \hfill{$\rhd$ update structural parameter estimates} \\
        $\bm{P}^{-1}=\left( \mathring{\bm{\Omega}}^\top \bm{C}^\top\bm{C} \mathring{\bm{\Omega}} \right)^{-1}$; \\
        ${\varrho} = {n-1} - \bm{1}^\top\bm{C} \mathring{\bm{\Omega}} \bm{P}^{-1} \mathring{\bm{\Omega}}^\top  \bm{C}^\top\bm{1} $; \\
        $\hat{\bm{\xi}}_i(\text{binds}) = \bm{P}^{-1} \mathring{\bm{\Omega}}^\top \bm{C}^\top
            \left\{ \bm{I}_{n-1} - \frac{1}{{\varrho}} \bm{1} \bm{1}^\top
                    \left[ \bm{I}_{n-1} - \bm{C} \mathring{\bm{\Omega}} \bm{P}^{-1} \mathring{\bm{\Omega}}^\top \bm{C}^\top \right]
            \right\} \bm{S}\mathring{\bm{X}}(:,i)$; \\
        }
        $\hat{\eta}_i = \frac{1}{{\varrho}} \bm{1}^\top \left[ \bm{I}_{n-1} - \bm{C} \mathring{\bm{\Omega}} \bm{P}^{-1} \mathring{\bm{\Omega}}^\top \bm{C}^\top \right] \bm{S}\mathring{\bm{X}}(:,i)$;
            \hfill{$\rhd$ generate initial condition estimate} \\
    }

    \Output{sparse structural estimates $\hat{\bm{\Xi}}=\left[\hat{\bm{\xi}}_1, \cdots, \hat{\bm{\xi}}_d\right]$, initial condition estimates $\hat{\bm{\eta}}=\left[\hat{\eta}_1, \cdots, \hat{\eta}_d\right]^\top$}
    \caption{Pseudo-code for STLS algorithm}
    \label{alg:1}
\end{algorithm}

\section{Simulations}\label{sec:3}

In this section the finite sample performance of ISINDy is investigated via Monte Carlo simulations; three nonlinear ODEs, namely Logistic, Lokta-Volterra, and Lorenz dynamics, are used as examples. The underlying time series $\bm{x}(t)$ is approximated with Runge-Kutta numerically. For each component of state variable, the noise level is measured by the noise-variance ratio
\begin{equation}\label{eq:19}
    nvr =
        \sqrt{\frac{\mathrm{Var}\left(\bm{e}_i\right)}{\mathrm{Var}\left(\bm{x}_i\right)}}
            \times 100\%
        = \frac{\sigma_i}{\sqrt{\mathrm{Var}\left(\bm{x}_i\right)}}
        \times 100\%,~
        i= 1,\cdots,d
\end{equation}
where $\bm{e}_i=\left[ e_i(t_1),~ e_i(t_2),~ \cdots,~ e_i(t_n) \right]^\top$ and $\bm{x}_i=\left[ x_i(t_1),~ x_i(t_2),~ \cdots,~ x_i(t_n) \right]^\top$.

{The proposed ISINDy is compared to the benchmark methods: the standard SINDy \cite{brunton2016discovering} and its integral variant presented in \cite{schaeffer2017sparse} (which is referred as to InSINDy hereafter). Note that these two benchmark methods do not consider measurement noises of the observations and also cannot produce estimations of initial condition. }

\subsection{Logistic equation }\label{sec:3-1}

Consider the logistic growth equation
\[
    \frac{d}{dt}x(t)=1.6{x}(t)-[{x}(t)]^2
\]
with initial condition $x(0)=\eta=0.1$.

\begin{figure}[!ht]
    \centering
    \includegraphics[scale=0.45]{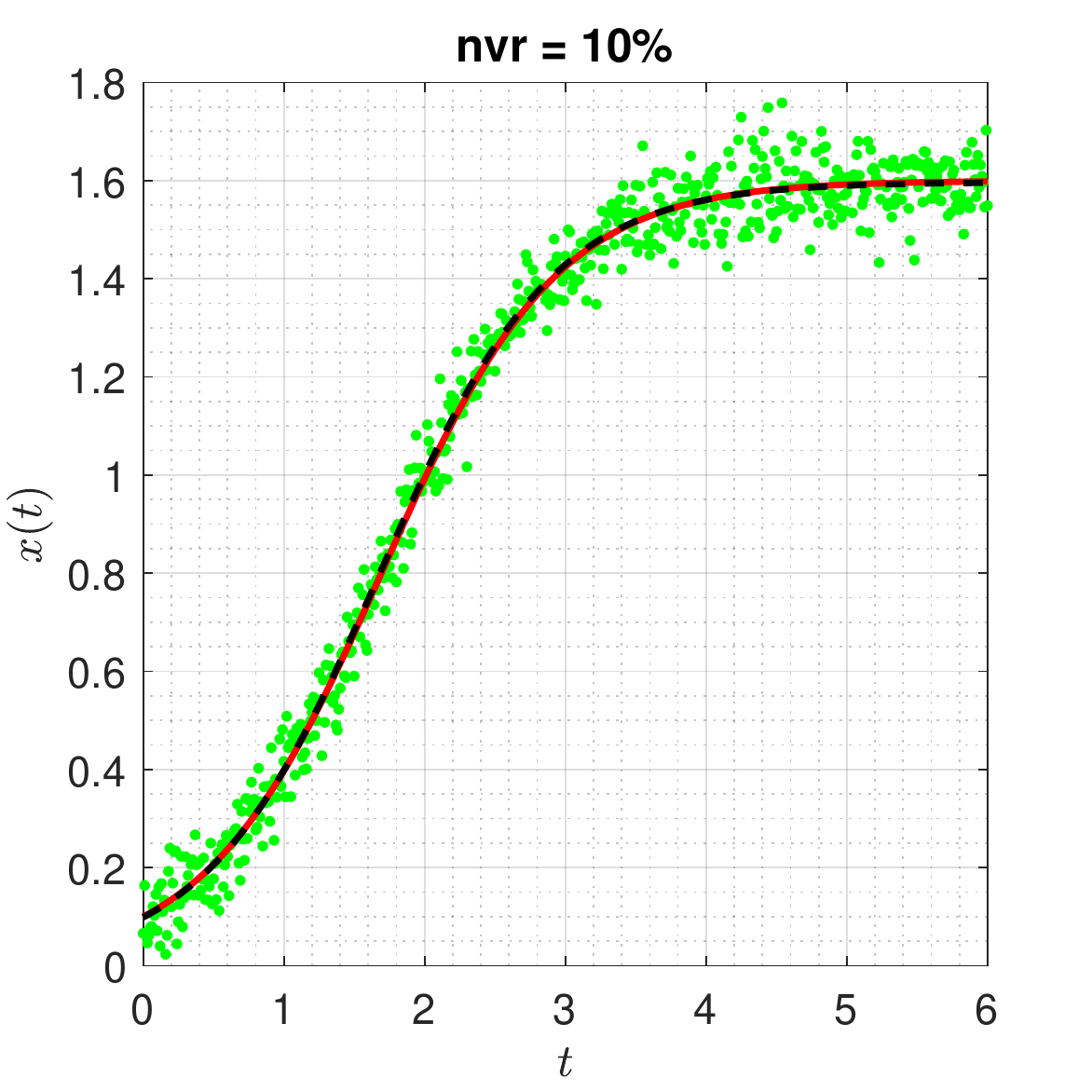}
    \includegraphics[scale=0.45]{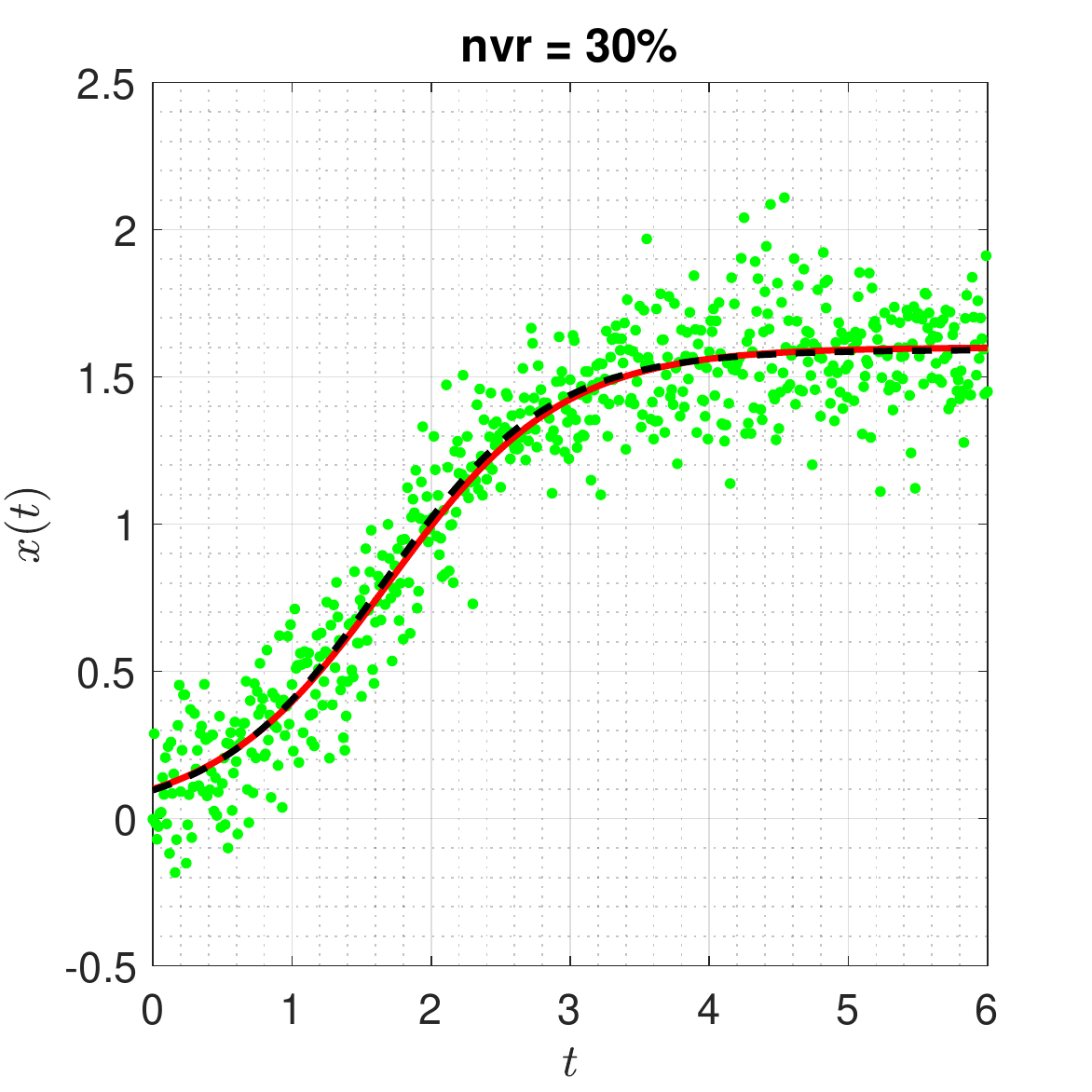}
    \includegraphics[scale=0.45]{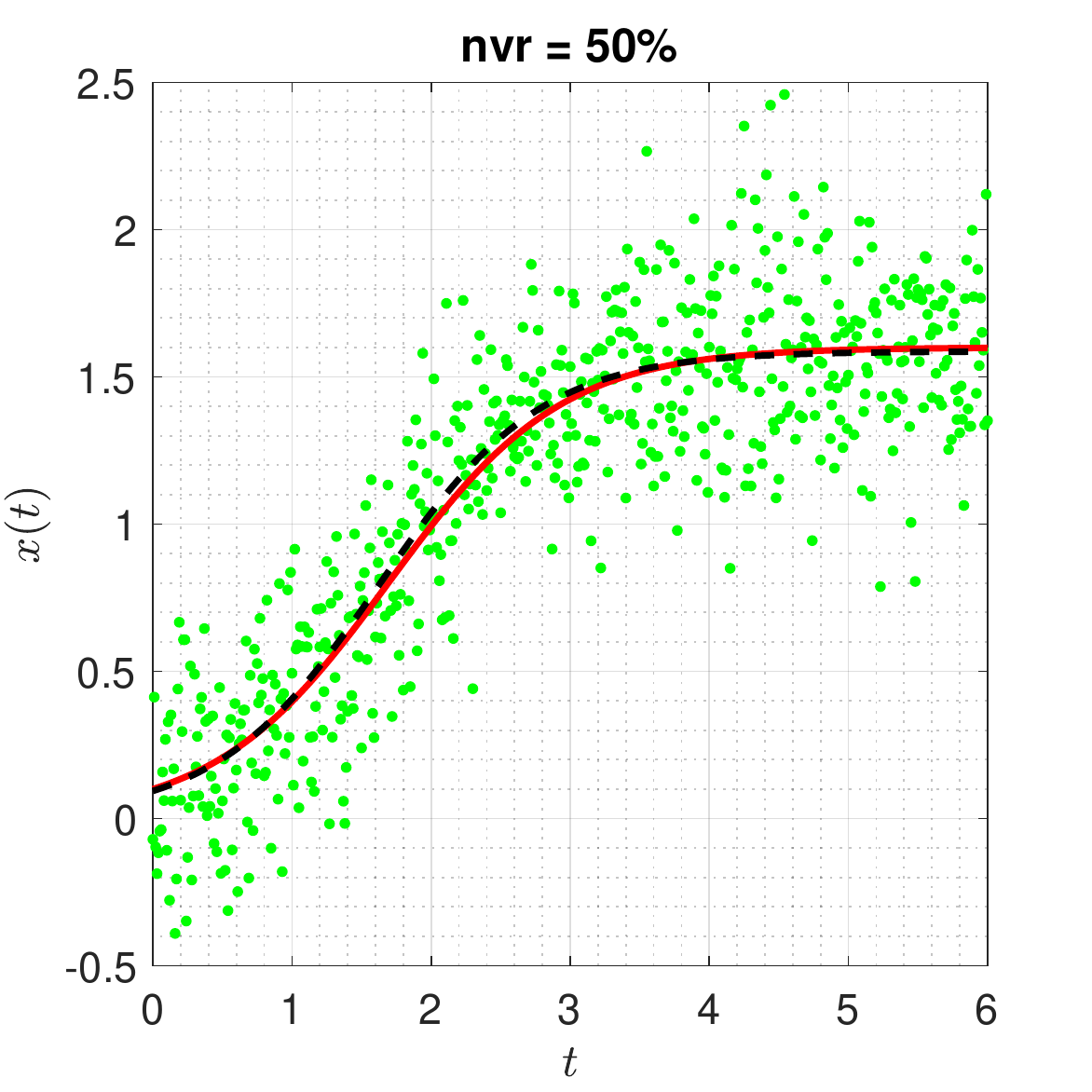} \\
    \caption{Noisy observations, true and identified trajectories of logistic equation. The noisy observations are shown in green dots; the true and identified trajectories are specified line style in red and dashed line style in black, respectively.}
    \label{fig:2}
\end{figure}

The trajectories are generated under the configuration: $t\in[0,6]$ with $h=0.01$, and the noise level $nvr\in[10\%, 30\%, 50\%]$. When implementing Algorithm \ref{alg:1}, the threshold value is set to $\lambda=0.1$. { In each noisy scenario, the estimates of sparse structural parameters and initial condition are summarised in Table \ref{tbl:1} (Appendix B), and the results show that only the proposed ISINDy identifies the correct active terms together with high-accuracy estimates. By substituting the ISINDy-based estimates and running the \textsf{ode45} routine, one has the identified trajectories in Figure \ref{fig:2}.
}

{
Figure \ref{fig:2} shows that the identified trajectories is very close to the true ones, even though the noise level is as high as $nvr=50\%$, indicating that ISINDy is robust to measurement noises. Moreover, both InSINDy and ISINDy can identify the correct active terms when the noise level is $nvr=10\%$, but ISINDy produces more accurate estimates and as the noise level increases, InSINDy fails to generate correct active terms. For example, when the noise level is $nvr=30\%$, the identified equations of ISINDy and InSINDy are
\[
    \frac{d}{dt}x(t)= 1.6609{x}(t) -1.0434[{x}(t)]^2
\]
and
\[
    \frac{d}{dt}x(t)= 2.3752{x}(t) -2.1302[{x}(t)]^2 + 0.4054[{x}(t)]^3,
\]
respectively. It is obvious that the InSINDy-based equation includes a redundant active term in the vector field.
}

\subsection{Lokta-Volterra system }\label{sec:3-2}

Consider the Lokta--Volterra system
\[
    \begin{cases}
        \frac{d}{dt}{x}_{1}(t) = \frac{2}{3}x_{1}(t)-\frac{4}{3}x_{1}(t)x_{2}(t) \\
        \frac{d}{dt}{x}_{2}(t) = -x_{2}(t)+x_{1}(t)x_2(t)
    \end{cases}
\]
with initial condition $[x_1(0), x_2(0)]^\top=[\eta_1, \eta_2]^\top=[1.8, 1.8]^\top$.

\begin{figure}[!ht]
    \centering
    \includegraphics[scale=0.4]{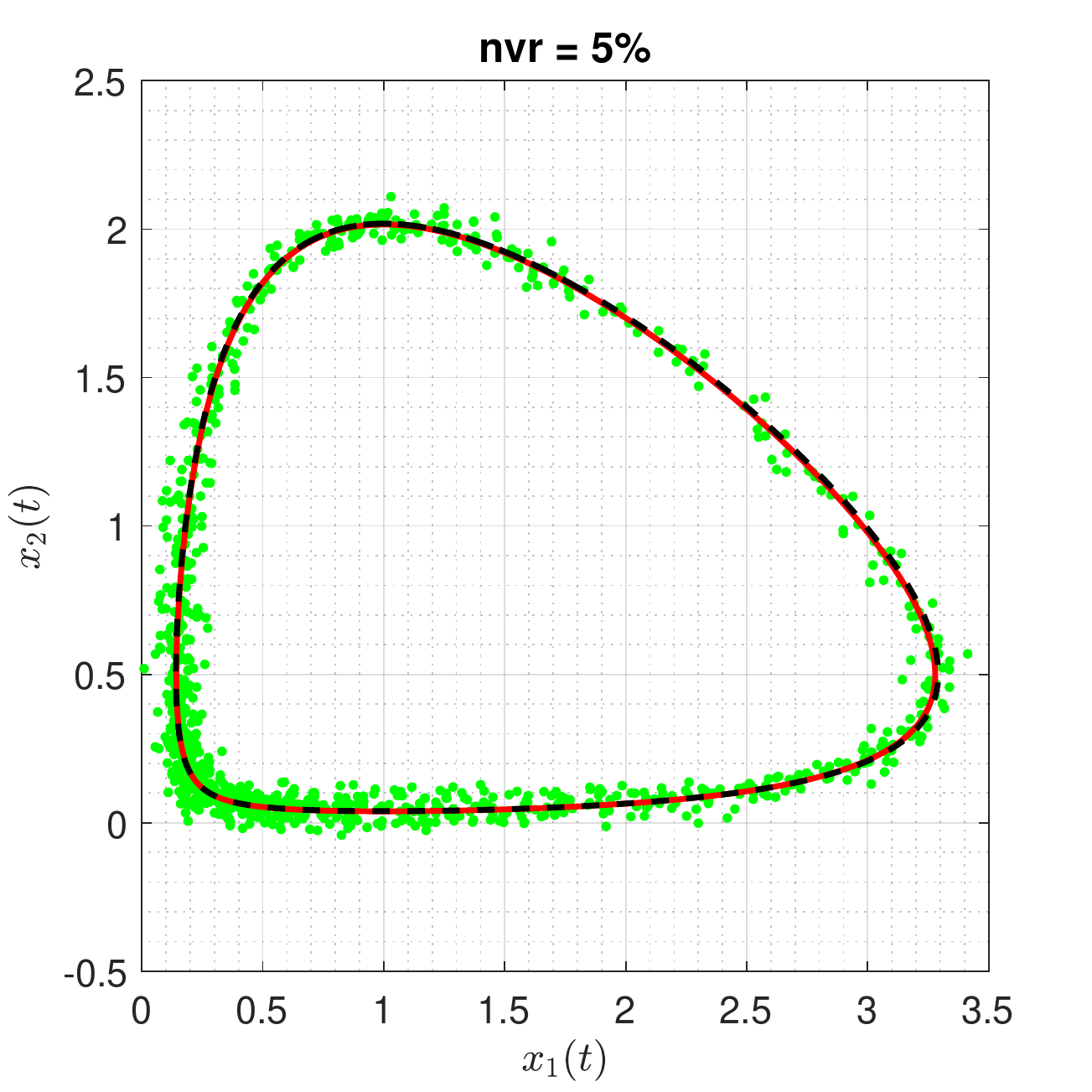}
    \includegraphics[scale=0.4]{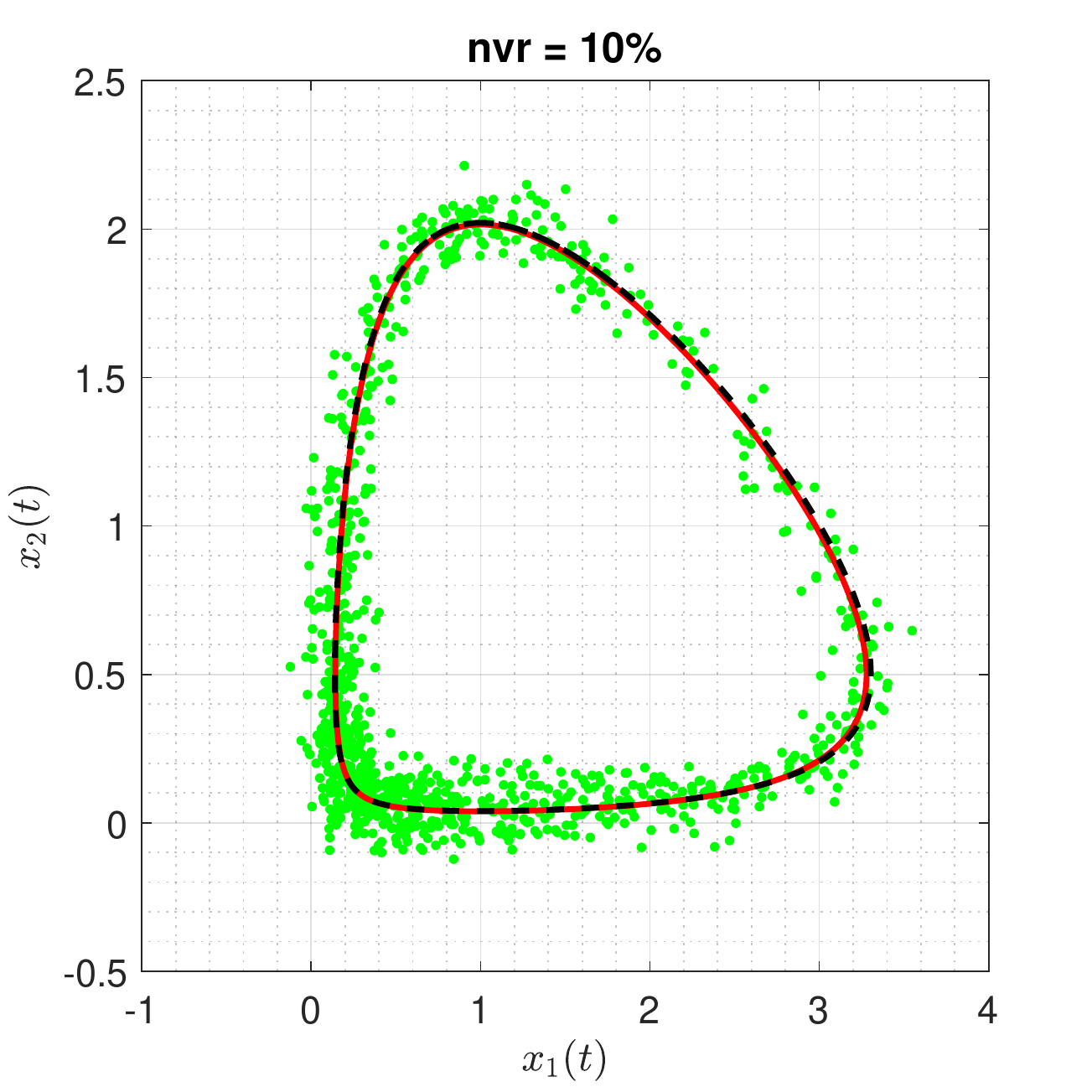}
    \includegraphics[scale=0.4]{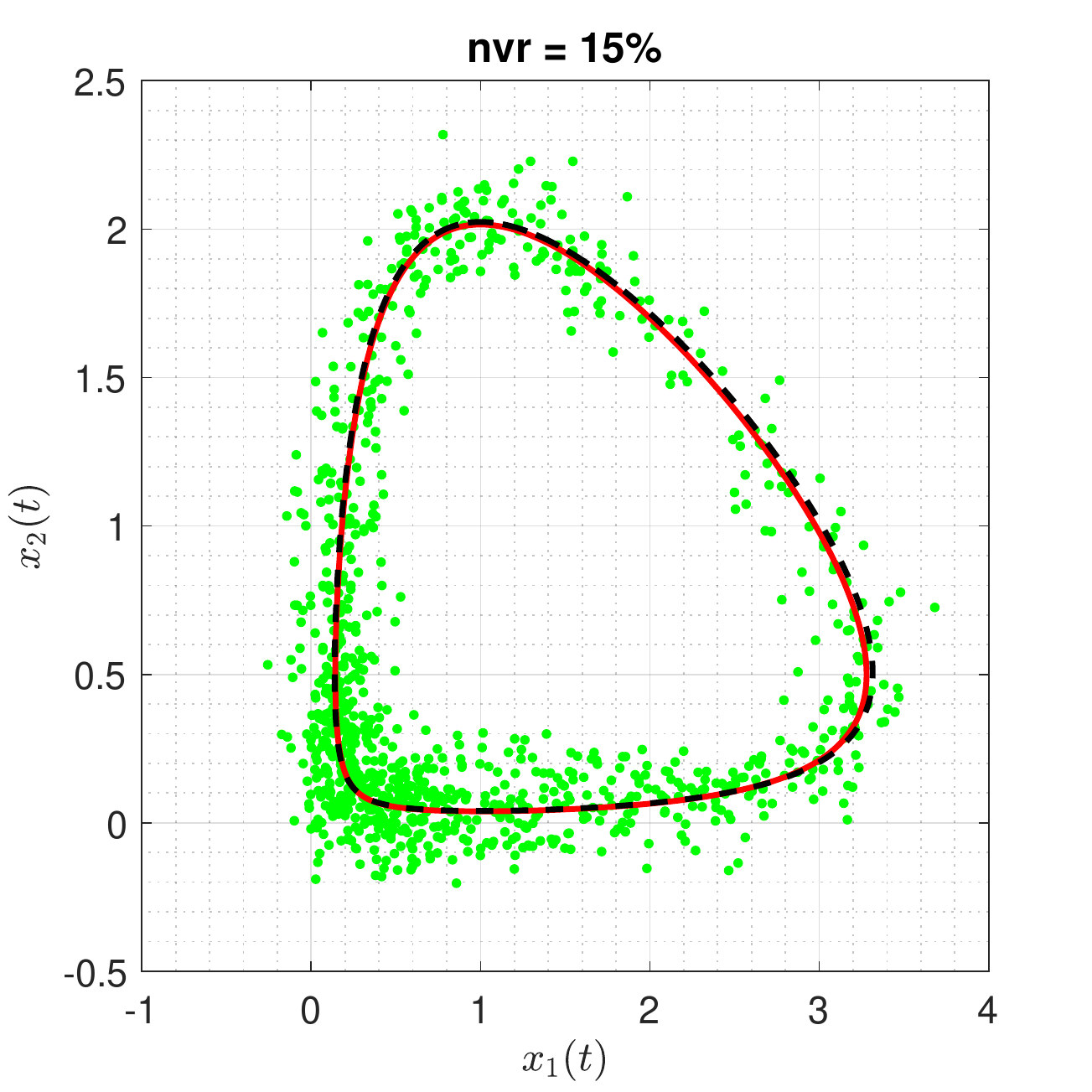} \\
    \caption{Noisy observations, true and identified trajectories of Lokta-Volterra system. The noisy observations are shown in green dots; the true and identified trajectories are specified line style in red and dashed line style in black, respectively.}
    \label{fig:3}
\end{figure}

The true trajectory are generated under the configuration: $t\in[0,10]$ with $h=0.01$, and the noise level $nvr\in[5\%, 10\%, 15\%]$. Setting the threshold value to $\lambda=0.3$, the estimates of sparse structural parameters and initial conditions are summarised in Table \ref{tbl:2} (Appendix B). {The results show that in all three noisy scenarios, only the proposed ISINDy identifies the correct active terms and generates high-accuracy parameter estimates.} By substituting these ISINDy-based estimates and running the \textsf{ode23} routine, one obtains the identified trajectories in Figure \ref{fig:3}.

\begin{figure}[!ht]
    \centering
    \includegraphics[scale=0.73]{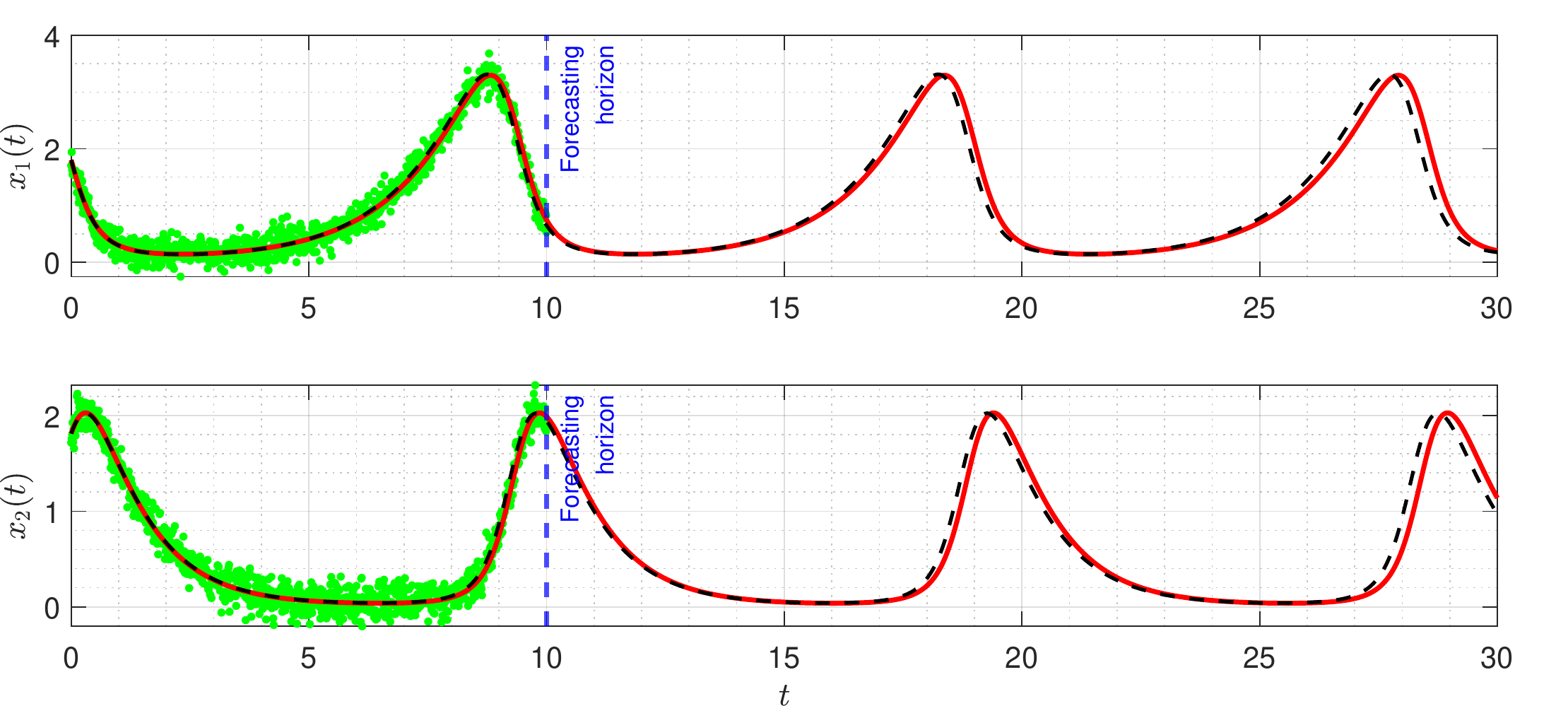} \\
    \caption{Noisy observations, true and identified trajectories of Lokta-Volterra system when noise level is $nvr=15\%$. The noisy observations are shown in green; the true and identified trajectories are specified line style in red and dashed line style in black, respectively.}
    \label{fig:3-1}
\end{figure}

Figure \ref{fig:3} shows that the proposed ISINDy captures the correct active terms in all three noisy scenarios, even though the noise level is up to $nvr=15\%$. {Increasing the noise level leads to larger errors in the estimation of structural parameters and initial condition, resulting in more biased trajectories.} Besides, to compare the identified and true dynamic systems in the presence of measurement noises on the observations, we simulate the true and identified trajectories at the time range $t\in[0,30]$. Figure \ref{fig:3-1} shows that the trajectory of the identified dynamical system almost coincides with the true one, indicating that ISINDy can also generate accurate long-term forecasts.

\subsection{Lorenz system }\label{sec:3-3}

Consider the chaotic Lorenz system
\[
    \begin{cases}
        \frac{d}{dt}{x}_{1}(t) = -10 x_1(t) +10 x_2(t)  \\
        \frac{d}{dt}{x}_{2}(t) = 28 x_1(t) - x_1(t) x_3(t) - x_2(t) \\
        \frac{d}{dt}{x}_{3}(t) = x_1(t) x_2(t)-\frac{8}{3} x_3(t) \\
    \end{cases}
\]
with initial condition $[x_1(0), x_2(0), x_3(0)]^\top=[\eta_1, \eta_2, \eta_3]^\top=[-5, 10, 30]^\top$.

\begin{figure}[!ht]
    \centering
    \includegraphics[scale=0.4]{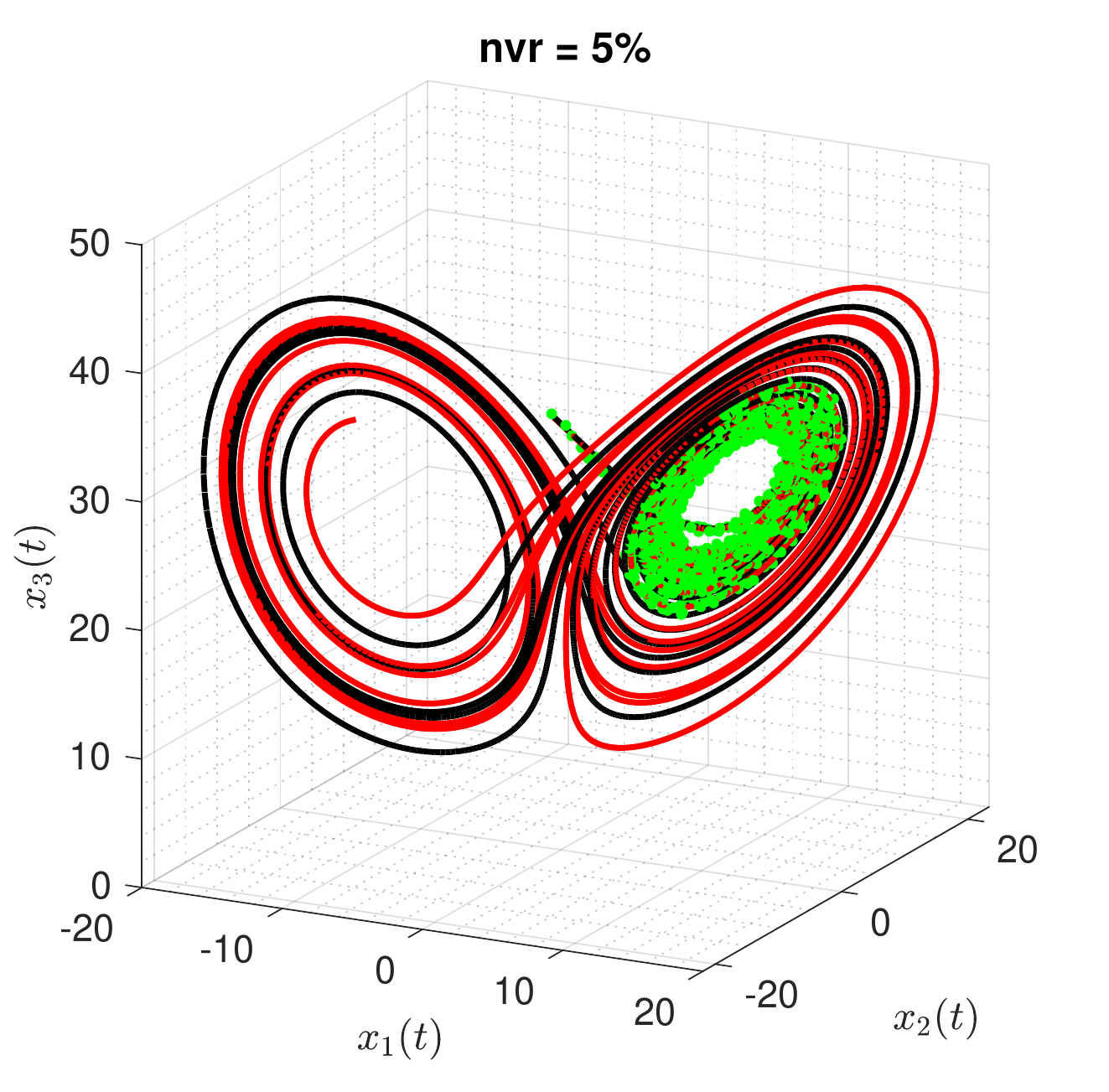}
    \includegraphics[scale=0.4]{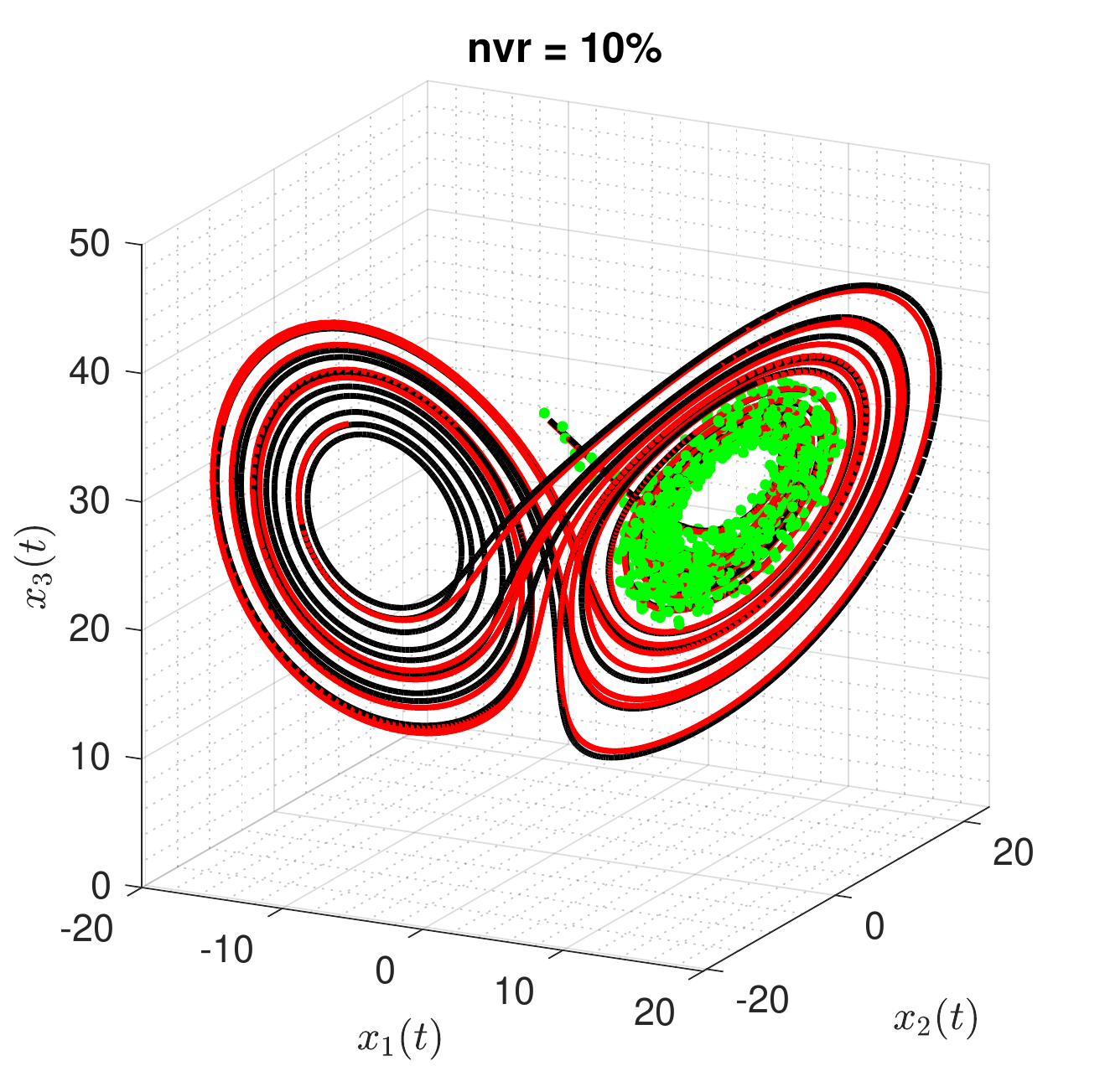}
    \includegraphics[scale=0.4]{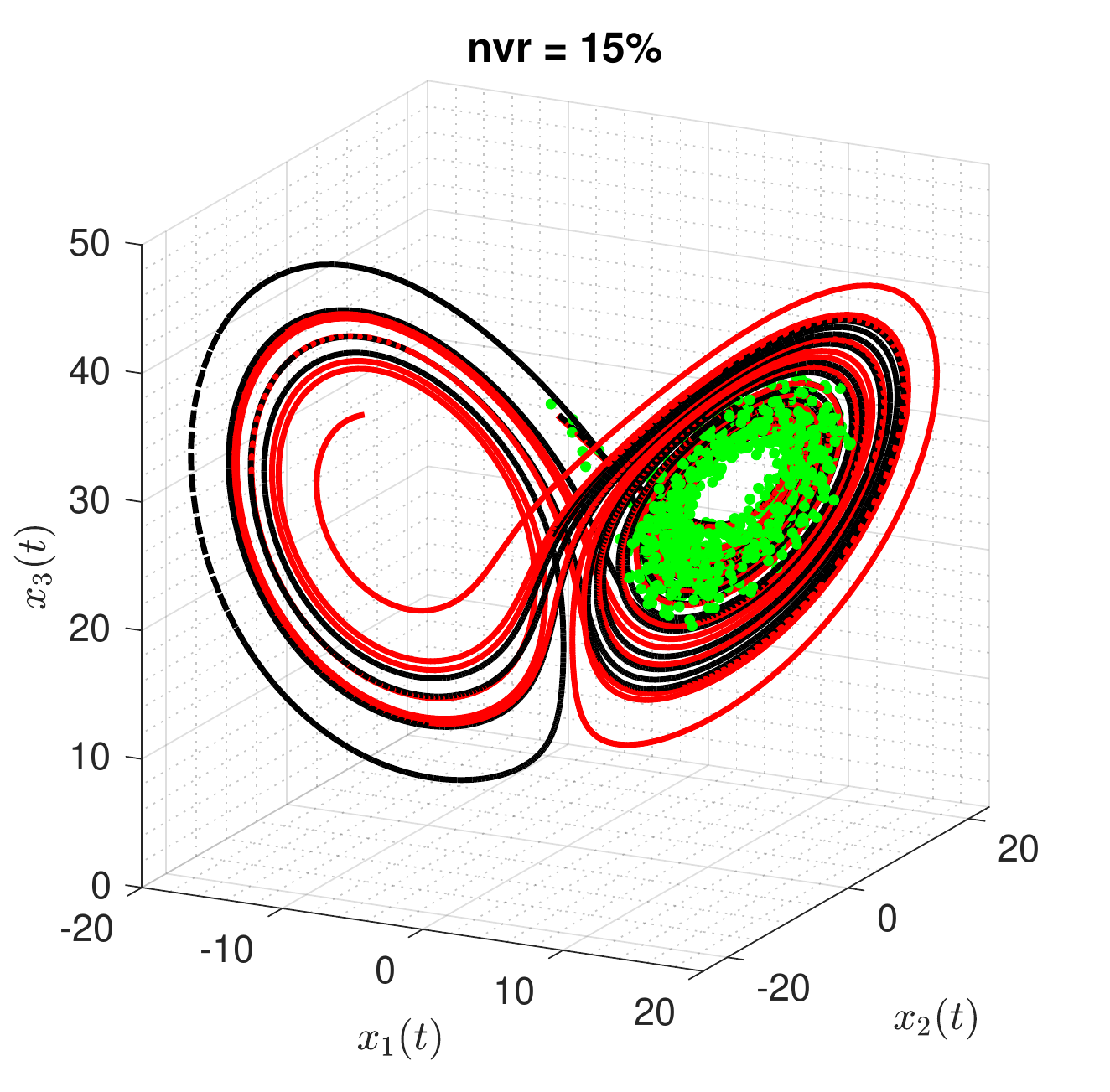} \\
    \caption{Noisy observations, true and identified trajectories of Lorenz system. The noisy observations are shown in green; the true and identified trajectories are specified line style in red and dashed line style in black, respectively.}
    \label{fig:4}
\end{figure}

The true trajectories are generated under the configuration: $t\in[0,5]$, $h=0.005$, and $nvr\in[5\%, 10\%, 15\%]$. The threshold value is set to $\lambda=0.8$, and the estimates of sparse structural parameters and initial conditions are summarised in Table \ref{tbl:3} (Appendix C). {The results show that SINDy, InSINDy, and ISINDy correctly identify the active terms in all three noisy scenarios for the first component, but SINDy and InSINDy fail to identify the correct active terms in the second and third components, letting alone the accurate estimates of parameters. Only ISINDy picks out the correct active terms and generates high-accuracy estimates in all three noisy scenarios.} By substituting the ISINDy-based estimates and running the \textsf{ode45} routine, one obtains the identified trajectories in Figure \ref{fig:4}.

\begin{figure}[!ht]
    \centering
    \includegraphics[scale=0.73]{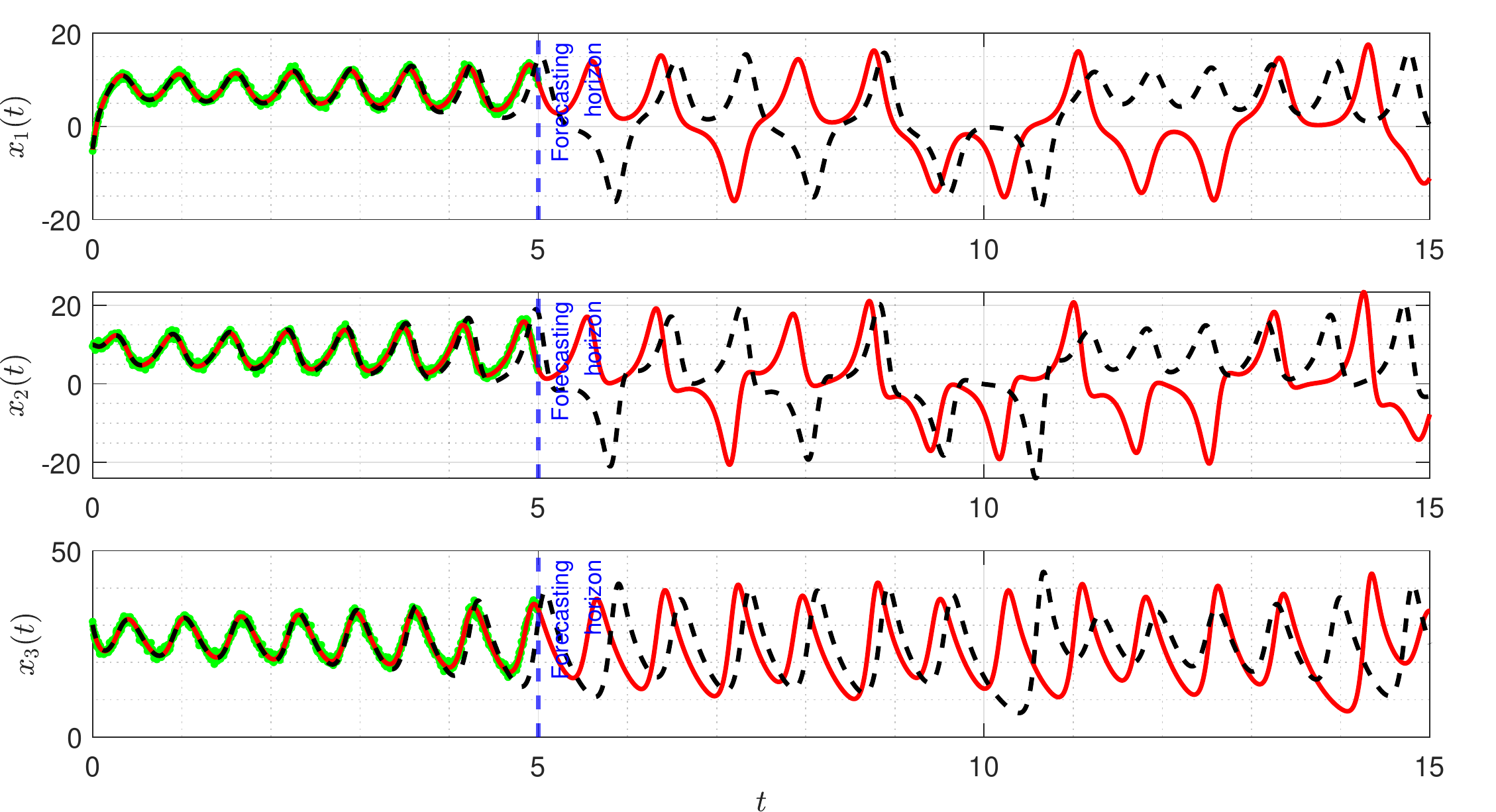} \\
    \caption{Noisy observations, true and identified trajectories of Lorenz system when noise level is $nvr=15\%$. The noisy observations are shown in green; the true and identified trajectories are specified line style in red and dashed line style in black, respectively.}
    \label{fig:5}
\end{figure}

\begin{figure}[!ht]
    \centering
    \includegraphics[scale=0.45]{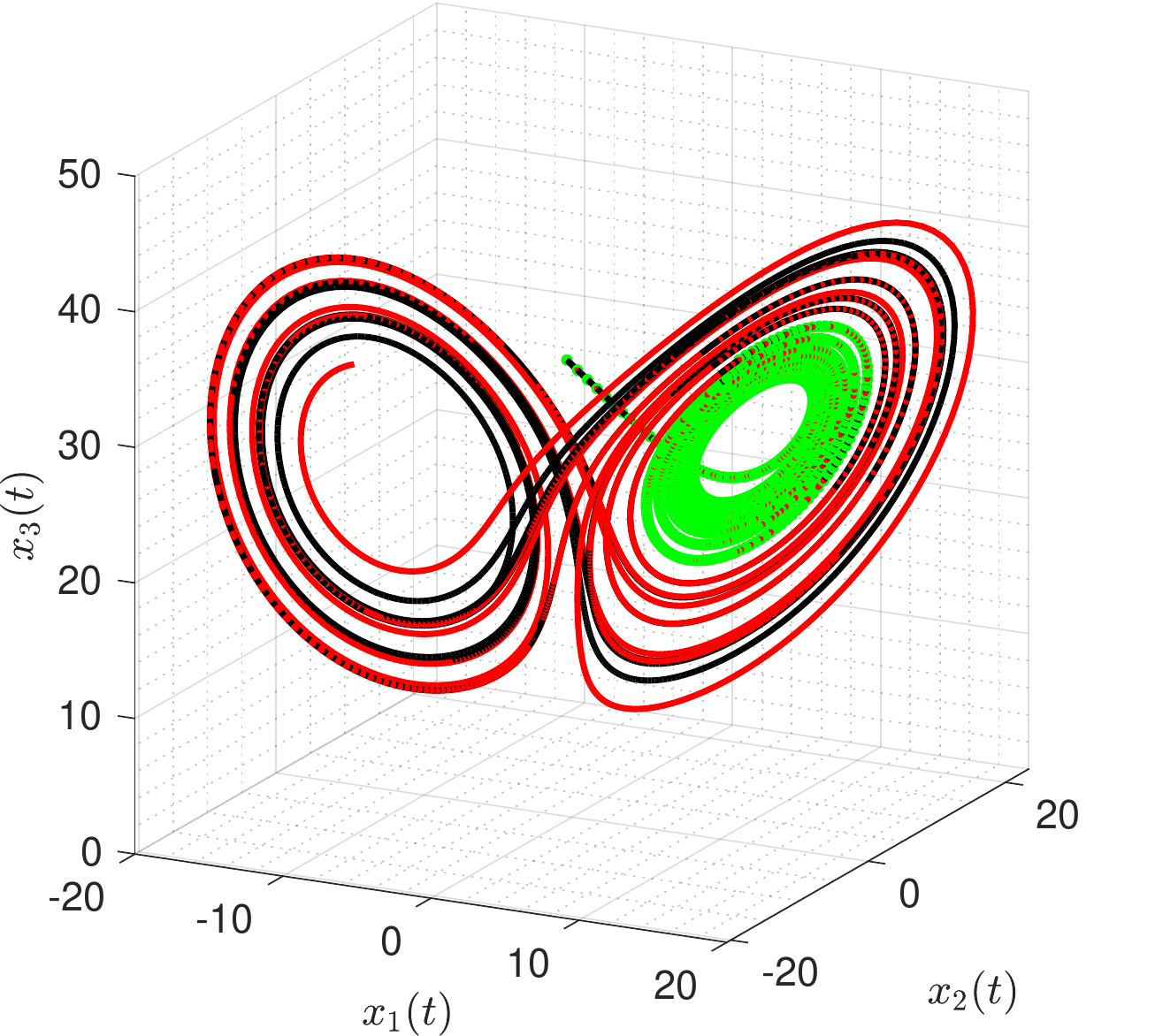}
    \includegraphics[scale=0.45]{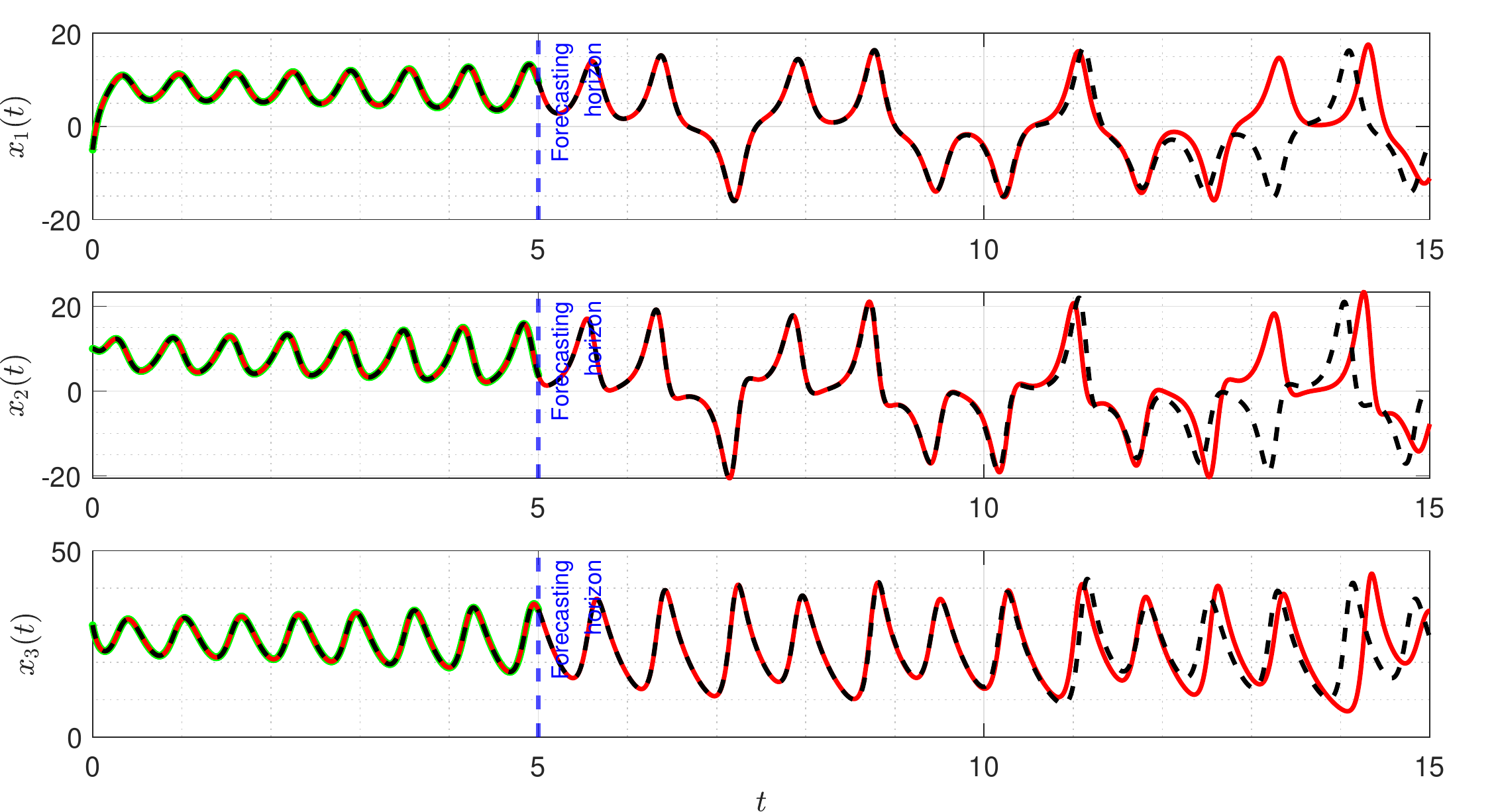} \\
    \caption{True and identified trajectories of Lorenz system when the observations are noise-free. The observations are shown in green; the true and identified trajectories are specified line style in red and dashed line style in black, respectively.}
    \label{fig:6}
\end{figure}

It is worth noting, however, that although the identified equation has correct model structure and also high-accuracy estimates, the identified trajectories may still stray off the true one quickly. For instance, in the $nvr=15\%$ scenario, the trajectories of three state variables are plotted in Figure \ref{fig:5}. It is clear that the identified trajectories almost coincides with the true ones for a short time from $t = 0$ to about $t = 4$ and then stray from about $t = 5$ on. The reason for this phenomenon is that Lorenz system has chaotic solutions for certain structural parameters and initial conditions. The identified trajectory is very sensitive to the estimates of structural parameters and initial condition. A small bias on the estimates results in large offsets in the trajectory. Specifically, assuming the observations are noise-free, i.e. $nvr=0$, the identified attractor associated with the trajectories of three states are shown in Figure \ref{fig:6}. The identified trajectories starts to fast stray off the true ones from about $t = 12$ on, indicating unreliable long-term forecasts even under a noise-free environment.

{
\subsection{Short discussions }

The results of three cases show that the proposed ISINDy outperforms the standard SINDy and its integral variant InSINDy, in both structure identification and parameter estimation terms. When observations are noise-free, i.e. $nvr=0\%$, all three methods can identify the correct model structure, but the parameter estimates obtained from ISINDy are much closer to their true values than those obtained form SINDy and InSINDy (see Tables \ref{tbl:1}--\ref{tbl:3} in Appendix B). To be specific, SINDy and ISINDy outperform InSINDy due to that the formers use a trapezoid formula rather than the Euler's formula, but this does not hold in noisy scenarios due to that the derivative approximation-based SINDy is sensitive to measurement noises and thus generates poor results, especially in the Lorenz system case (see Table \ref{tbl:3} in Appendix B).

It is worth noting, however, that the aforementioned simulations are conducted under two conditions: (1) the time range is fixed in a specified interval, and (2) the candidate features include all possible terms appearing in vector fields. Actually, these two factors have important effects on the performance of ISINDy.
\begin{enumerate}
    \item  [(i)]
    For factor (1), taking the Lorenz system as an example, the fixed time interval is set to $[0,3]$, $[3,7]$, and $[7,12]$, respectively. When the observations are noise-free, ISINDy identifies the correct model structure in all three time interval settings (see Figure \ref{fig:c} in Appendix C). But, when the observations are contaminated by noises, the noise robustness of ISINDy in each interval is different. Figure \ref{fig:d} in Appendix C shows that ISINDy generates wrong model structure when $nvr=2\%$, indicating that in $[0,3]$ the noise tolerance is less than 2\%. Similarly, in $[3, 7]$ the noise tolerance is less than 5\% and in $[7, 12]$ the noise tolerance is less than 10\%. An alternative strategy to improve robustness is guaranteeing the quality of smoothing, which is inferred from the noise-free observations cases.

    \item [(ii)]
    For factor (2), a simple ODE $\frac{d}{dt} x = -\sin(x )$ is employed to explore the effect of candidate features on identification accuracy. Note that the Taylor expansion of the vector field is
    \[
        g(x) = -\sin(x) = - x + \frac{1}{3!} x^3 - \frac{1}{5!} x^5 + \mathcal{O}(x^7)
                        = - x + 0.16667 x^3 - 0.0083 x^5 + \mathcal{O}(x^7).
    \]
    On this test problem, ISINDy with different candidate features is investigated, including polynomial, trigonometric, and a combination of both, shown in Table \ref{tbl:4} (Appendix D). In the case of polynomial features, the degree of polynomial is set to 3 and 5 in turn, and in both settings ISINDy identifies the correct active terms together with high-accuracy parameters in the Taylor expansion of $-\sin(x)$. In both of the other cases of trigonometric features or a combination of polynomial and trigonometric features, the correct active term $\sin(x)$ is identified.

\end{enumerate}

Furthermore, the time interval and the design of candidate features always interact. Taking the simple ODE $\frac{d}{dt}x=-\sin(x)$ as an example, in the time range of $t\in [0, 5]$, the resulting trajectories lie in the range of $x(t)\in [-1.0, 1.0]$. This guarantees that the Taylor series expansion can approximate the vector filed, thereby generating the excellent dynamic reconstruction, as shown in Figure \ref{fig:e} (Appendix D).
}

\section{Conclusion }\label{sec:4}

In this work, an integral SINDy approach that leverages penalized spline smoothing is proposed for simultaneously (1) denoising time-series observations, (2) generating accurate estimates of structural parameters and initial condition, and (3) identifying parsimonious systems of differential equations responsible for generating time series. Using the integral formulation, noisy time-series observations are handled in a stable and robust way; even for high noise levels, the algorithm can produce correct model structure, high-accuracy estimates of structural parameters and initial condition, showing the improvement of proposed approach with respect to robustness to noise with less training data in comparison with the previous integral variant of SINDy \cite{schaeffer2017sparse}. Simulations show that this approach is accurate and robust to noise when dynamic behavior is included in the observations, suggesting its ability to extract interpretable equations from data.

{ The current research assumes that all terms appeared in the vector field are included in the candidate features.} We must acknowledge that the proposed approach has several limitations, e.g., (1) time-consuming of penalized spline smoothing process, especially for large-sample observations, and (2) incorrect identification of vector fields if the candidate features are improperly designed. {In practice, how to design an appropriate library of candidate features is still a challenging problem \cite{delahunt2022toolkit}, although physical intuition or domain knowledge may be leveraged in many systems.}
Besides, note that this work does not relate to the identifiability of structural parameters and initial condition \cite{piasaccomani2003parameter,miao2011identifiability,peng2011conditions}, and the uncertainty quantification of model structure inference \cite{degennaro2019model}, which are the centre topic of our future theoretical research.

\section*{Data availability }

All scripts to reproduce the results are available at \url{https://github.com/weibl9/integralSINDy}.

%

\section*{Declaration of competing interest }

The author declares no competing interests.

\bibliography{SparseDynamics}

%

{

\section*{Appendix A. The cubic B-spline basis functions }

A cubic B-spline function is formed by jointing piece-wise polynomials of degree 3 with at most $\mathcal{C}^2$ continuity at the knots. The basis functions between two knots are defined by the Cox-de Boor recursion formula \cite{patrikalakis2002representation}:
\begin{align*}
    & b_{i,0}(t) = \begin{cases}
            1, & t \in    [\tau_i, \tau_{i+1} ) \\
            0, & t \notin [\tau_i, \tau_{i+1} )
        \end{cases} \\
    & b_{i,p}(t) = \frac{t-\tau_i}{\tau_{i+p}-\tau_i}b_{i,p-1}(t) + \frac{\tau_{i+p+1}-t}{\tau_{i+p+1}-\tau_{i+1} }b_{i+1,p-1}(t)
\end{align*}
where $\tau_i$ denotes knots and $p=3$ denotes degree of polynomial.

For example, let the knots be $\{0, 0.1, \cdots, 1.0\}$. Then, 13 basis functions are defined over the interval $[0,1]$, and at each knot the
polynomial values and their first two derivatives match, as show in Figure \ref{fig:a-1}.

\begin{figure}[!ht]
    \centering
    \includegraphics[scale=1.0]{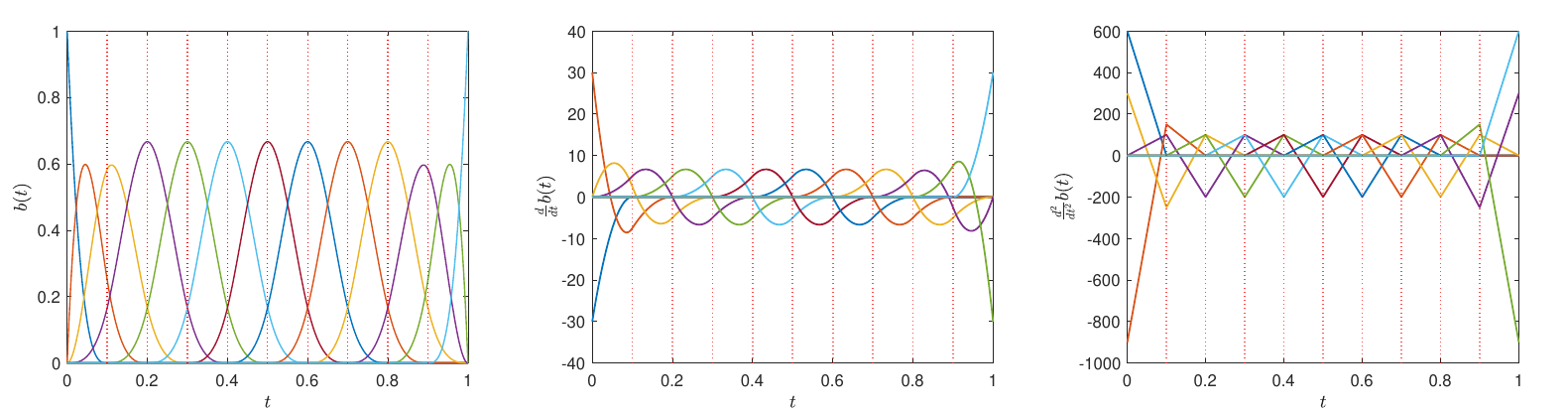}
    \caption{Left panel shows the basis functions, and middle (resp. right) panel shows the first-order (resp. second-order) derivatives of basis functions. }
    \label{fig:a-1}
\end{figure}

\cleardoublepage
\section*{Appendix B. The estimates of structural parameters and initial condition under varying noise levels }

\begin{table}[!ht]
    \centering
    \footnotesize
    \caption{Sparse structural parameters and initial condition of logistic equation under varying noise levels. }
    \label{tbl:1}
    \setlength{\tabcolsep}{1.5mm}
    \begin{tabular}{l rrr rrr rrr rrr} \toprule
    & \multicolumn{3}{c}{$nvr = 0\%$} & \multicolumn{3}{c}{$nvr = 10\%$} & \multicolumn{3}{c}{$nvr = 30\%$} & \multicolumn{3}{c}{$nvr = 50\%$ } \\
        \cmidrule(r){2-4} \cmidrule(r){5-7} \cmidrule(r){8-10} \cmidrule(r){11-13}
    Terms  & SINDy & InSINDy & ISINDy & SINDy & InSINDy & ISINDy & SINDy & InSINDy & ISINDy & SINDy & InSINDy & ISINDy \\ \hline
$\eta$ & --- & --- & 0.1000 & --- & --- & 0.0985 & --- & --- & 0.0962 & --- & --- & 0.0935 \\
$x$ & 1.6000 & 1.5975 & 1.6000 & 1.8604 & 1.6831 & 1.6203 & 3.7477 & 2.3752 & 1.6609 & 5.1682 & 3.4485 & 1.7047 \\
$x^2$ & $-$1.0000 & $-$0.9980 & $-$1.0000 & $-$1.4848 & $-$1.0551 & $-$1.0145 & $-$5.0504 & $-$2.1302 & $-$1.0434 & $-$7.4255 & $-$3.5429 & $-$1.0744 \\
$x^3$ & 0 & 0 & 0 & 0.2016 & 0 & 0 & 1.6942 & 0.4054 & 0 & 2.5584 & 0.8550 & 0 \\ \bottomrule
    \end{tabular}
\end{table}

\begin{table}[!ht]
    \centering
    \footnotesize
    \caption{Sparse structural parameters and initial condition of Lokta-Volterra equation under varying noise levels.  }
    \label{tbl:2}
    \setlength{\tabcolsep}{1.5mm}
    \begin{tabular}{ll rrr rrr rrr rrr} \toprule
& & \multicolumn{3}{c}{$nvr = 0\%$} & \multicolumn{3}{c}{$nvr = 5\%$} & \multicolumn{3}{c}{$nvr = 10\%$} & \multicolumn{3}{c}{$nvr = 15\%$ } \\
    \cmidrule(r){3-5} \cmidrule(r){6-8} \cmidrule(r){9-11} \cmidrule(r){12-14}
  & Terms  & SINDy & InSINDy & ISINDy & SINDy & InSINDy & ISINDy & SINDy & InSINDy & ISINDy & SINDy & InSINDy & ISINDy \\ \hline
First & $\eta_1$ & --- & --- & 1.8000 & --- & --- & 1.8017 & --- & --- & 1.8030 & --- & --- & 1.8032 \\
Comp. & $x_1$ & 0.6667 & 0.6692 & 0.6667 & 0.6664 & 0.6698 & 0.6720 & 0.6589 & 0.6702 & 0.6774 & 0.6425 & 0.5575 & 0.6825 \\
& $x_1^2$ & 0 & 0 & 0 & 0 & 0 & 0 & 0 & 0 & 0 & 0 & 0 & 0 \\
& $x_1^3$ & 0 & 0 & 0 & 0 & 0 & 0 & 0 & 0 & 0 & 0 & 0 & 0 \\
& $x_2$ & 0 & 0 & 0 & 0 & 0 & 0 & 0 & 0 & 0 & 0 & $-$0.7683 & 0 \\
& $x_1x_2$ & $-$1.3334 & $-$1.3430 & $-$1.3333 & $-$1.3292 & $-$1.3277 & $-$1.3388 & $-$1.3114 & $-$1.3122 & $-$1.3444 & $-$1.2771 & 2.7419 & $-$1.3488 \\
& $x_1^2x_2$ & 0 & 0 & 0 & 0 & 0 & 0 & 0 & 0 & 0 & 0 & $-$0.9956 & 0 \\
& $x_2^2$ & 0 & 0 & 0 & 0 & 0 & 0 & 0 & 0 & 0 & 0 & 0.7534 & 0 \\
& $x_1x_2^2$ & 0 & 0 & 0 & 0 & 0 & 0 & 0 & 0 & 0 & 0 & $-$0.9186 & 0 \\
& $x_2^3$ & 0 & 0 & 0 & 0 & 0 & 0 & 0 & 0 & 0 & 0 & $-$0.4481 & 0 \\ \hline
Second & $\eta_2$ & --- & --- & 1.8000 & --- & --- & 1.8037 & --- & --- & 1.8077 & --- & --- & 1.8121 \\
Comp. & $x_1$ & 0 & 0 & 0 & 0 & 0 & 0 & 0 & 0 & 0 & 0 & 0 & 0 \\
& $x_1^2$ & 0 & 0 & 0 & 0 & 0 & 0 & 0 & 0 & 0 & 0 & 0 & 0 \\
& $x_1^3$ & 0 & 0 & 0 & 0 & 0 & 0 & 0 & 0 & 0 & 0 & 0 & 0 \\
& $x_2$ & $-$0.9999 & $-$0.9968 & $-$1.0000 & $-$1.1301 & $-$1.0956 & $-$0.9954 & $-$0.9973 & $-$0.9677 & $-$0.9908 & $-$0.8007 & $-$1.0828 & $-$0.9863 \\
& $x_1x_2$ & 0.9990 & 0.9978 & 1.0000 & 1.7761 & 0 & 0.9956 & 1.8408 & 0.9851 & 0.9911 & 1.9741 & 0 & 0.9865 \\
& $x_1^2x_2$ & 0 & 0 & 0 & 0 & 0 & 0 & 0 & 0 & 0 & 0 & 0 & 0 \\
& $x_2^2$ & 0 & 0 & 0 & 0 & 0 & 0 & 0 & 0 & 0 & $-$0.5601 & 0 & 0 \\
& $x_1x_2^2$ & 0 & 0 & 0 & 0 & 0.7589 & 0 & 0 & 0 & 0 & $-$0.4198 & 0.7636 & 0 \\
& $x_2^3$ & 0 & 0 & 0 & 0 & 0 & 0 & 0 & 0 & 0 & 0 & 0 & 0 \\ \bottomrule
\end{tabular}
\end{table}

\begin{table}[!ht]
    \centering
    \caption{Sparse structural parameters and initial condition of Lorenz system with varying noise levels. Note that the rows with all entries equal to 0 are omitted here. }
    \label{tbl:3}
    \footnotesize
    \setlength{\tabcolsep}{1mm}
    \begin{tabular}{ll rrr rrr rrr rrr} \toprule
    & & \multicolumn{3}{c}{$nvr = 0\%$} & \multicolumn{3}{c}{$nvr = 5\%$} & \multicolumn{3}{c}{$nvr = 10\%$} & \multicolumn{3}{c}{$nvr = 15\%$ } \\
    \cmidrule(r){3-5} \cmidrule(r){6-8} \cmidrule(r){9-11} \cmidrule(r){12-14}
    & Terms  & SINDy & InSINDy & ISINDy & SINDy & InSINDy & ISINDy & SINDy & InSINDy & ISINDy & SINDy & InSINDy & ISINDy \\ \hline
    First & $\eta_1$ & --- & --- & $-$5.0046 & --- & --- & $-$4.9798 & --- & --- & $-$4.9215 & --- & --- & $-$4.8467 \\
    Comp. & $x_1$ & $-$10.0013 & $-$10.2411 & $-$10.0013 & $-$10.0200 & $-$10.3430 & $-$10.0180 & $-$9.9638 & $-$10.4420 & $-$10.0100 & $-$9.8377 & $-$10.5390 & $-$9.9915 \\
    & $x_2$ & 10.0012 & 10.2435 & 10.0013 & 10.0000 & 10.3270 & 9.9987 & 9.9286 & 10.4070 & 9.9728 & 9.7906 & 10.4840 & 9.9353 \\
    & $\vdots$ & $\vdots$ & $\vdots$ & $\vdots$ & $\vdots$ & $\vdots$ & $\vdots$ & $\vdots$ & $\vdots$ & $\vdots$ & $\vdots$ & $\vdots$ & $\vdots$ \\
    & $x_3^3$       & 0 & 0 & 0 & 0 & 0 & 0 & 0 & 0 & 0 & 0 & 0 & 0 \\ \hline
    Second & $\eta_2$ & --- & --- & 9.9983 & --- & --- & 9.9822 & --- & --- & 9.9451 & --- & --- & 9.8771 \\
    Comp. & $x_1$ & 27.9523 & 27.9406 & 28.0061 & 59.9790 & 25.8860 & 28.0770 & 52.5500 & 26.0150 & 28.0720 & 66.2880 & 26.1250 & 27.9900 \\
    & $x_1^2$       & 0 & 0 & 0 & 5.1065 & 0 & 0 & 1.8500 & 0 & 0 & $-$2.4958 & 0 & 0 \\
    & $x_2$ & $-$0.9820 & $-$0.9683 & $-$1.0002 & 7.9855 & 0 & $-$0.9592 & 47.4620 & 0 & $-$0.8955 & 68.9500 & 0 & $-$0.8065 \\
    & $x_1x_2$      & 0 & 0 & 0 & $-$3.5886 & 0 & 0 & 5.3775 & 0 & 0 & 12.8470 & 0 & 0 \\
    & $x_2^2$ & 0 & 0 & 0 & $-$1.0011 & 0 & 0 & $-$5.1637 & 0 & 0 & $-$7.4029 & 0 & 0 \\
    & $x_3$ & 0 & 0 & 0 & $-$12.214 & 0 & 0 & $-$27.9080 & 0 & 0 & $-$44.8650 & 0 & 0 \\
    & $x_1x_3$ & $-$0.9970 & $-$0.9990 & $-$1.0002 & $-$5.5552 & $-$0.96061 & $-$1.0047 & $-$7.0064 & $-$0.9654 & $-$1.0073 & $-$9.1726 & $-$0.9695 & $-$1.0079 \\
    & $x_2x_3$ & 0 & 0 & 0 & 1.0570 & 0 & 0 & $-$1.4457 & 0 & 0 & $-$3.2327 & 0 & 0 \\
    & $x_3^2$ & 0 & 0 & 0 & 0.9825 & 0 & 0 & 2.4698 & 0 & 0 & 4.0385 & 0 & 0 \\
    & $\vdots$ & $\vdots$ & $\vdots$ & $\vdots$ & $\vdots$ & $\vdots$ & $\vdots$ & $\vdots$ & $\vdots$ & $\vdots$ & $\vdots$ & $\vdots$ & $\vdots$ \\
    & $x_3^3$       & 0 & 0 & 0 & 0 & 0 & 0 & 0 & 0 & 0 & 0 & 0 & 0 \\ \hline
    Third & $\eta_3$ & --- & --- & 30.0043 & --- & --- & 30.0420 & --- & --- & 30.0700 & --- & --- & 30.0950 \\
    Comp. & $x_1$ & 0 & 0 & 0 & 50.8140 & 0 & 0 & 96.1630 & 1.8294 & 0 & 130.3000 & $-$10.6670 & 0 \\
    & $x_1^2$ & 0 & 0 & 0 & 5.3069 & 0 & 0 & 12.2670 & 0 & 0 & 17.3530 & 0 & 0 \\
    & $x_2$ & 0 & 0 & 0 & $-$58.8090 & 0 & 0 & $-$121.9800 & 0 & 0 & $-$157.9900 & 0 & 0 \\
    & $x_1x_2$ & 0.9949 & 1.0088 & 1.0002 & $-$8.3655 & 1.0142 & 0.9964 & $-$21.4980 & 0.9183 & 0.9923 & $-$29.0380 & 1.1835 & 0.9879 \\
    & $x_2^2$ & 0 & 0 & 0 & 4.5377 & 0 & 0 & 9.6943 & 0 & 0 & 12.2090 & 0 & 0 \\
    & $x_3$ & $-$2.6710 & $-$2.6934 & $-$2.6672 & 6.1171 & $-$2.7172 & $-$2.6625 & 15.8890 & $-$3.0193 & $-$2.6572 & 17.9310 & 0 & $-$2.6507 \\
    & $x_1x_3$ & 0 & 0 & 0 & $-$3.5181 & 0 & 0 & $-$6.5122 & 0 & 0 & $-$9.5474 & 0 & 0 \\
    & $x_2x_3$ & 0 & 0 & 0 & 4.0361 & 0 & 0 & 9.3839 & 0 & 0 & 12.7760 & 0 & 0 \\
    & $x_1x_2x_3$ & 0 & 0 & 0 & 0 & 0 & 0 & 0 & 0 & 0 & 1.1170 & 0 & 0 \\
    & $x_3^2$ & 0 & 0 & 0 & 0 & 0 & 0 & $-$1.5798 & 0 & 0 & $-$1.7850 & 0 & 0 \\
& $\vdots$ & $\vdots$ & $\vdots$ & $\vdots$ & $\vdots$ & $\vdots$ & $\vdots$ & $\vdots$ & $\vdots$ & $\vdots$ & $\vdots$ & $\vdots$ & $\vdots$ \\
    & $x_3^3$       & 0 & 0 & 0 & 0 & 0 & 0 & 0 & 0 & 0 & 0 & 0 & 0  \\ \bottomrule
    \end{tabular}
\end{table}

\cleardoublepage
\section*{Appendix C. The identified trajectories of Lorenz system under varying time intervals }

\begin{figure}[!ht]
    \centering
    \includegraphics[scale=0.95]{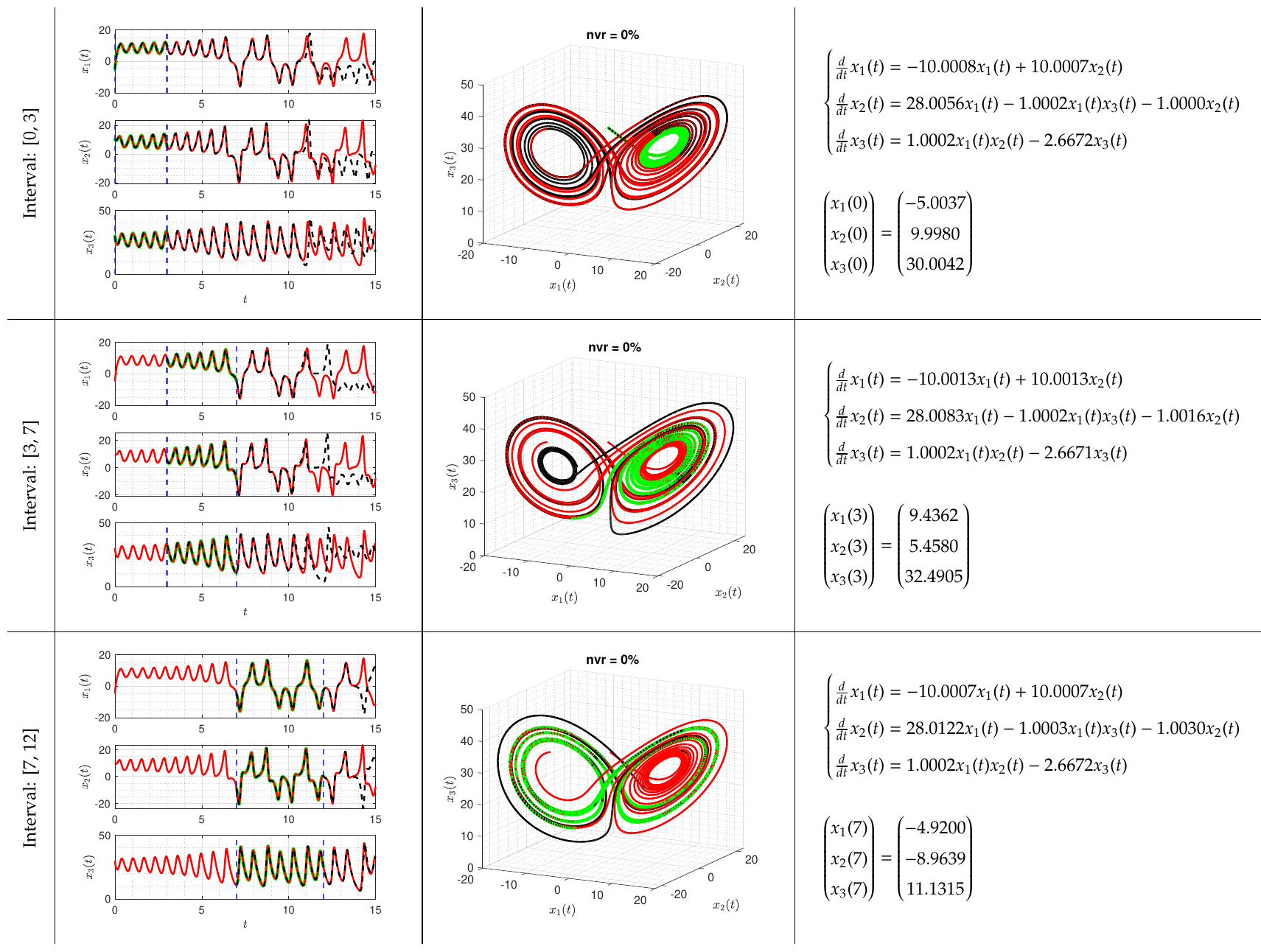}
    \caption{Noise-free condition: true and identified trajectories of Lorenz system under varying time intervals. The observations are shown in green; the true and identified trajectories are specified line style in red and dashed line style in black, respectively. }
    \label{fig:c}
\end{figure}

\begin{figure}[!ht]
    \centering
    \includegraphics[scale=0.85]{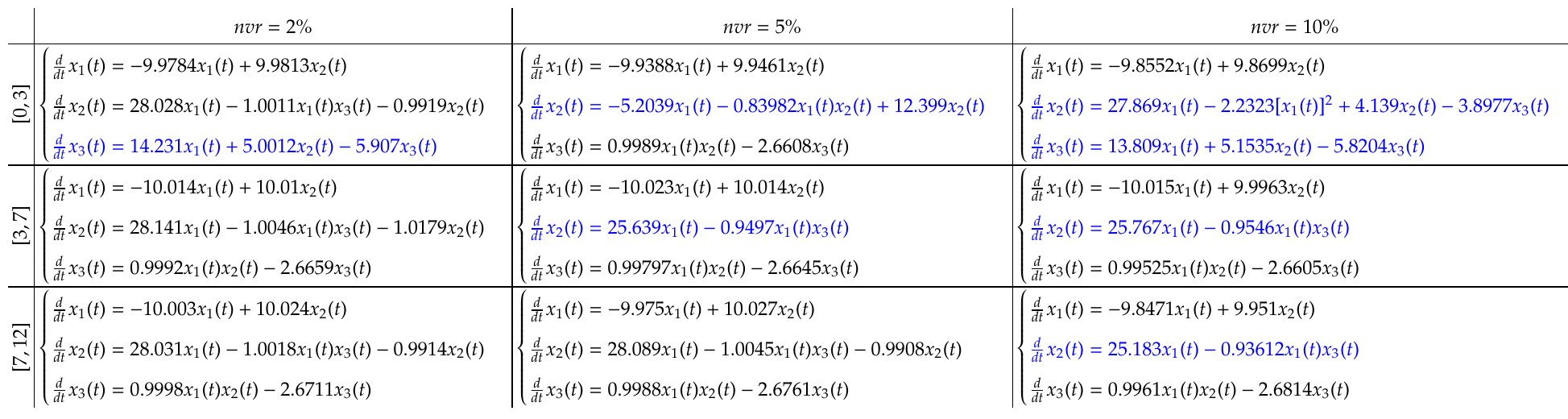}
    \caption{Noisy condition: identified equations of Lorenz system under varying time intervals. The misspecified components are in blue. }
    \label{fig:d}
\end{figure}

\cleardoublepage
\section*{Appendix D. The effect of candidate features on structure identification }

\begin{table}[!ht]
    \centering
    \caption{Sparse structural parameters and initial condition of Lorenz system with varying noise levels. }
    \label{tbl:4}
    \footnotesize
    \setlength{\tabcolsep}{2.6mm}
    \begin{tabular}{lrrrrrrrrrr} \toprule
\shortstack{initial\\ condition} & $-$1.00 & $-$0.80 & $-$0.60 & $-$0.40 & $-$0.20 & 0.20 & 0.40 & 0.60 & 0.80 & 1.00 \\ \hline
\multicolumn{11}{l}{polynomial features: degree 3} \\
$\eta$ & $-$0.9997 & $-$0.7999 & $-$0.6000 & $-$0.4000 & $-$0.2000 & 0.2000 & 0.4000 & 0.6000 & 0.7999 & 0.9997 \\
$x$ & $-$0.9996 & $-$0.9998 & $-$0.9999 & $-$1.0000 & $-$1.0000 & $-$1.0000 & $-$1.0000 & $-$0.9999 & $-$0.9998 & $-$0.9996 \\
$x^2$ & 0 & 0 & 0 & 0 & 0 & 0 & 0 & 0 & 0 & 0 \\
$x^3$ & 0.1612 & 0.1632 & 0.1647 & 0.1658 & 0.1664 & 0.1664 & 0.1658 & 0.1647 & 0.1632 & 0.1612 \\ \hline
\multicolumn{11}{l}{polynomial features: degree 5} \\
$\eta$ & $-$1.0000 & $-$0.8000 & $-$0.6000 & $-$0.4000 & $-$0.2000 & 0.2000 & 0.4000 & 0.6000 & 0.8000 & 1.0000 \\
$x$ & $-$1.0000 & $-$1.0000 & $-$1.0000 & $-$1.0000 & $-$1.0000 & $-$1.0000 & $-$1.0000 & $-$1.0000 & $-$1.0000 & $-$1.0000 \\
$x^2$ & 0 & 0 & 0 & 0 & 0 & 0 & 0 & 0 & 0 & 0 \\
$x^3$ & 0.1666 & 0.1666 & 0.1666 & 0.1666 & 0.1666 & 0.1666 & 0.1666 & 0.1666 & 0.1666 & 0.1666 \\
$x^4$ & 0 & 0 & 0 & 0 & 0 & 0 & 0 & 0 & 0 & 0 \\
$x^5$ & $-$0.0081 & $-$0.0082 & $-$0.0082 & $-$0.0083 & $-$0.0083 & $-$0.0083 & $-$0.0083 & $-$0.0082 & $-$0.0082 & $-$0.0081 \\ \hline
\multicolumn{11}{l}{trigonometric features} \\
$\eta$ & $-$1.0000 & $-$0.8000 & $-$0.6000 & $-$0.4000 & $-$0.2000 & 0.2000 & 0.4000 & 0.6000 & 0.8000 & 1.0000 \\
$\sin(x)$ & $-$1.0000 & $-$1.0000 & $-$1.0000 & $-$1.0000 & $-$1.0000 & $-$1.0000 & $-$1.0000 & $-$1.0000 & $-$1.0000 & $-$1.0000 \\
$\cos(x)$ & 0 & 0 & 0 & 0 & 0 & 0 & 0 & 0 & 0 & 0 \\
$\sin(2x)$ & 0 & 0 & 0 & 0 & 0 & 0 & 0 & 0 & 0 & 0 \\
$\cos(2x)$ & 0 & 0 & 0 & 0 & 0 & 0 & 0 & 0 & 0 & 0 \\ \hline
\multicolumn{11}{l}{combination of polynomial and trigonometric features} \\
$\eta$ & $-$1.0000 & $-$0.8000 & $-$0.6000 & $-$0.4000 & $-$0.2000 & 0.2000 & 0.4000 & 0.6000 & 0.8000 & 1.0000 \\
$x$ & 0 & 0 & 0 & 0 & 0 & 0 & 0 & 0 & 0 & 0 \\
$x^2$ & 0 & 0 & 0 & 0 & 0 & 0 & 0 & 0 & 0 & 0 \\
$x^3$ & 0 & 0 & 0 & 0 & 0 & 0 & 0 & 0 & 0 & 0 \\
$\sin(x)$ & $-$1.0000 & $-$1.0000 & $-$1.0000 & $-$1.0000 & $-$1.0000 & $-$1.0000 & $-$1.0000 & $-$1.0000 & $-$1.0000 & $-$1.0000 \\
$\cos(x)$ & 0 & 0 & 0 & 0 & 0 & 0 & 0 & 0 & 0 & 0 \\
$\sin(2x)$ & 0 & 0 & 0 & 0 & 0 & 0 & 0 & 0 & 0 & 0 \\
$\cos(2x)$ & 0 & 0 & 0 & 0 & 0 & 0 & 0 & 0 & 0 & 0 \\ \bottomrule
    \end{tabular}
\end{table}

\begin{figure}[!ht]
    \centering
    \includegraphics[scale = 1.0]{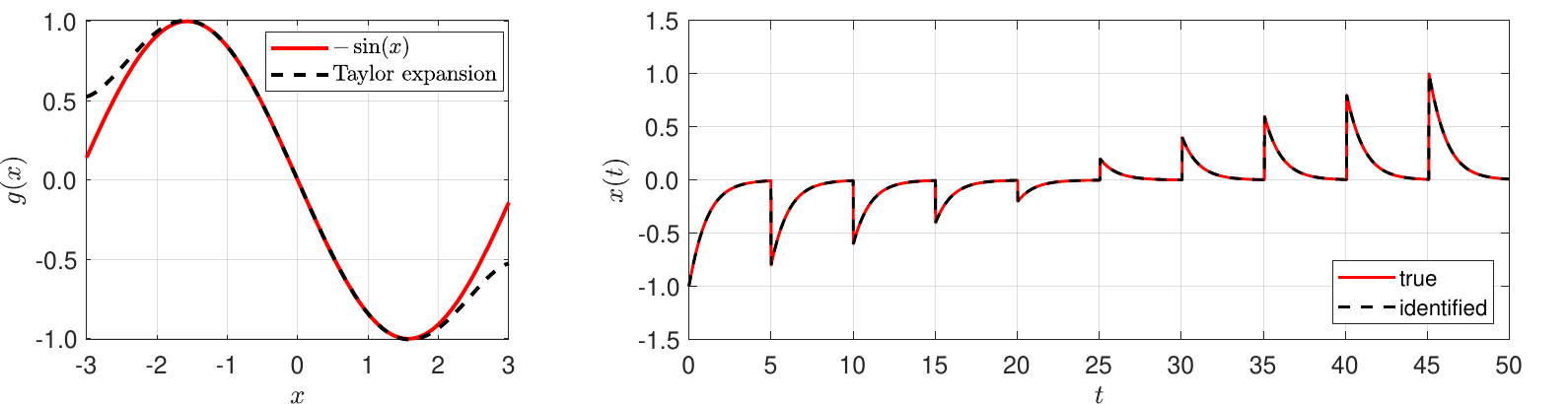}
    \caption{Left panel shows the Taylor series expansion of $-\sin(x)$ up to the fifth order in $[-\pi,\pi]$. Right panel shows the sparse dynamics reconstruction for a sequence of initial condition initialized every 5 time units. Initial conditions are chosen from $-$1.0 to 1.0 in increments of 0.2 (excluding 0.0). }
    \label{fig:e}
\end{figure}

}

\end{document}